\documentclass{article}
\usepackage[centertags]{amsmath}

\begin{document}

\title{ Average values of functionals and concentration without measure}

\author{ Cheng-shi Liu\\Department of Mathematics\\Northeast Petroleum University\\Daqing 163318, China
\\Email:chengshiliu-68@126.com}

\maketitle

Although there doesn't exist the Lebesgue measure in the ball $M$ of $C[0,1]$ with $p-$norm, the average values (expectation) $EY$ and variance $DY$ of some functionals $Y$ on $M$ can still be defined through the procedure of limitation from finite dimension to infinite dimension. In particular, the probability densities of coordinates of points in the ball $M$ exist and are derived out even though the density of points in $M$ doesn't exist. These densities include high order normal distribution, high order exponent distribution. This also can be considered as the geometrical origins of these probability distributions. Further, the exact values (which is represented in terms of finite dimensional integral) of a kind of infinite-dimensional functional integrals are obtained, and specially the variance $DY$ is proven to be zero, and then the nonlinear exchange formulas of average values of  functionals are also given. Instead of measure, the variance is used to measure the deviation of functional from its average value. $DY=0$ means that a functional takes its average on a ball with probability 1 by using the language of probability theory, and this is just the concentration  without measure. In addition, we prove that the average value depends on the discretization.


\section{Introduction}

In complexity science and  statistical physics in special, we often need to deal with high dimensional data and a large number of free degrees. Sometimes these data and free degrees can be considered as infinite-dimensional variables, and some physical quantities can be represented by infinite dimensional integrals.  Therefore, we need study the infinite-dimensional integrals.  The computations of integrals on functions with infinite number of
variables is still an important and interesting topic in quantum and statistical physics or even in finance(see, for example, [1-3]).  In 1920's, in the works of  G\^{a}teaux and  L\'{e}vy [4], the infinite-dimensional integrals had been considered and computed from the view of the point of probability theory. Further, Wiener integral became an important tool in stochastic processes theory[5-16]. Today, there exist a large number of papers devoting to the computations and applications of functional integrals. In particular, recently, some new algorithms such as  multilevel and changing dimension algorithms or dimension-wise quadrature methods, are proposed to approximate such integrals efficiently[17-29]. From the popular viewpoint, the foundation of functional integration such as Wiener's integral should be obtained on the theory of measure. However, Feynman's path integral is still lacking of a satisfied measure theory. On the other hand, there are some interesting problems in infinite dimensional space so that we have to consider infinite dimensional integrations under the condition of nonexistence of Lebesque measure, while other measure such as Gauss measure is not suitable to our aims. For example, if we randomly  take a continuous function $x(t)\in M$ where $M=\{x|a\leq x(t)\leq b, x(t)\in C[0,1]\}$, what is the average value of area $Y=\int_0^1x(t)\mathrm{d}t$? Here the randomness means that we take the points in $M$ by equal possibility, and hence we need a Lebesque measure on $M$. But it is well-known that there doesn't exist the Lebesque measure on $M$ at all,
While, conceptually, this problem is rather natural. Without measure theory, we can also use a limit procedure to give a rigorous definition of the average value (see below definition 1), and  the average value of $Y$ is easily  solved. In the paper,  we consider a more complex case in which $M$ is taken as the ball in $C[0,1]$ with $p-$norm, that is, $M=\{x|\int_0^1|x^p(t)|\mathrm{d}t\leq R, x(t)\in C[0,1]\}$.

Although there doesn't exist the Lebesque measure and then the density of points in $M$ doesn't exist, we show that the probability densities of coordinates of points in the ball $M$ do exist and are derived out with forms of high order normal and exponent distributions. Further,  we define and compute the exact average values (which are represented by the finite dimensional integrals)(expectation) and variances of some functionals. If we formally consider these functionals as the infinite-dimensional
random variables, the considered infinite-dimensional integral is just the expectation of the infinite-dimensional random variable. We show that the variances are zeros to prove that these functionals satisfy the property of the complete concentration of "measure". This is because if the measure exists, $DY=0$ means that a functional takes its average value on an infinite-dimensional ball with probability 1. In our cases of no Lebesque measure, we use $DY=0$ to replace the complete concentration of measure under the probability meaning.  Corresponding, we give the nonlinear exchange formulas for averages of functionals. The usual concentration of measure is described by some inequalities such as L\'{e}vy lemma[30-32], which is different with the complete concentration of "measure" which is shown by variance being zero.

Abstractly, a functional
$f(x)$ is a function of $x$ where $x$ is an element in an
infinite-dimensional space such as $C[0,1]$.  In
general, there are two basic ways to construct functionals. One
method is to use the values of $x(t)$ on some points
$t_1,\cdots,t_m$ such that
\begin{equation*}
f(x)=g(x(t_1),\cdots,x(t_m)),
\end{equation*}
where $g$ is a usual function in $R^n$. Essentially, such
functionals all are finite dimensional functions. Another method is to use
integral of $x(t)$ on some sets such that
\begin{equation*}
Y=f(x)=\int_{I_1}\cdots\int_{I_m}g(x(t_1),\cdots,x(t_m))\mathrm{d}t_1\cdots\mathrm{d}t_m,
\end{equation*}
where $I_1,\cdots,I_m$ are subsets of the interval $[0,1]$.  Such
functionals all are really infinite-dimensional functions.
Therefore, there are two kinds of basic elements $x(t_i)$ and
$\int_{I_i}g(x(t))\mathrm{d}t$ such that many interesting functionals
can be constructed in terms of them by addition, subtraction, multiplication, division, composition and limitation.

 For the first kind of functionals, the functional integral is just the usual finite-dimensional
integral. Thus  we only  consider the second
kind of functionals $f(x)$. If the domain of the functional $f$ is $M$, the
integral of $f$ on $M$ can be formally written as
\begin{equation}
\int_{M} f(x)D(x),
\end{equation}
where $D(x)$ represents formally the differential
element of the volume of $M$. But, in general, under the  meaning of Lebesque's measure, the volume $\int_{M}D(x)$ of $M$ is zero or infinity, and the infinite-dimensional integral $\int_{M}
f(x)D(x)$ is also respectively zero or infinity. However, the average value of functional $f$
on $M$,
\begin{equation}
Ef=\frac{\int_{M} f(x)D(x)}{\int_{M}D(x)}
\end{equation}
perhaps exists and is finite or infinite in general. Firstly, we  need a reasonable
 definition of the average value of functional. Since there doesn't exist the Lebesque measure in infinite-dimensional space in general, our approach is to use a limit procedure to define the average value of functionals.
 For example, similar to G\^{a}teaux and  L\'{e}vy  (see, [4]), we give the following definition.

\textbf{ Definition 1}: For $M=\{x|a\leq x(t)\leq b, x(t)\in C[0,1]\}$ and $Y=f(x)=\int_0^1g(x(t))\mathrm{d}t$ where $g(x)$ is a continuous function on $[a,b]$,
 we can define the average value of $f$ as
\begin{equation}
Ef=\lim_{n\rightarrow \infty}\frac{\int_a^b\cdots\int_a^b \frac{1}{n}\sum_{k=1}^ng(x_k)\mathrm{d}x_1\cdots\mathrm{d}x_n}{\int_a^b\cdots\int_a^b \mathrm{d}x_1\cdots\mathrm{d}x_n}
\end{equation}
where $x_k=x(\frac{k}{n})$. If the limitation exists and is finite or infinite, we call it the average value of functional $f$ on $M$.

\textbf{Remark 1}. Since $g(x)$ is continuous on $[a,b]$, for any $x_k=x(t_k)$ where $t_k\in[\frac{k}{n},\frac{k+1}{n})$, the above limitation is independent to the choice of $t_k$. Of important is that $Y$ and $M$ must take the same $t_k$ and $x_k$. From theorem 7 in section 4, we can see that the average value depends on the discretization! In addition, for $m$ variables function $g(x_1,\cdots,x_m)$, the average value of the functional $Y=f(x)=\int_0^1\cdots \int_0^1g(x(t_1),\cdots,x(t_m))\mathrm{d}t_1\cdots\mathrm{d}t_m$ can be defined similarly. More general, if $a=a(t), b=b(t)$ are two fixed continuous functions, we can give corresponding definition of the average value of $f$ for
$M=\{x|a(t)\leq x(t)\leq b(t), x(t)\in C[0,1]\}$,
\begin{equation*}
Ef=\lim_{n\rightarrow \infty}\frac{\int_{a_1}^{b_1}\cdots\int_{a_n}^{b_n} \frac{1}{n}\sum_{k=1}^ng(x_k)\mathrm{d}x_1\cdots\mathrm{d}x_n}{\int_{a_1}^{b_1}\cdots\int_{a_n}^{b_n} \mathrm{d}x_1\cdots\mathrm{d}x_n},
\end{equation*}
where $a_i=a(t_i), b_i=b(t_i), x_k=x(t_k)$.

Here, we must point out that the equiprobability (or equal possibility) hypothesis is implicated in the definition, that is, the points in $M$ are taken by equiprobability. This is a natural assumption. For example, for the aforementioned problem,  if we take randomly a continuous function $x(t)\in M =\{x|0\leq x(t)\leq 1\}, t\in[0,1]$, what is the average value of its area $Y=\int_0^1x(t)\mathrm{d}t$? In the problem, we have implicitly supposed that we take the function $x(t)$ in $M$ by equiprobability. However, this is just an intuitive and formal explanation in infinite-dimensional cases because in general there exists no  Lebesque measure on $M$ as the probability measure to give the meaning of equiprobability ([16]). But in finite dimensional cases this assumption has strict mathematical foundation since there exists the Lebesque's measure as the corresponding probability measure such that we can talk about reasonably equiprobability. In the whole paper, when we say equiprobability, it is just the meaning here. Below we give such definition.

\textbf{Definition 2}. Let $M$ be a bounded set in the space of all continuous functions on $[0,1]$  with some norm.
If for any finite-dimensional subset $M_0$ of $M$, the points in $M_0$ are taken by equiprobability, we say that the points in $M$ are taken by equiprobability (or equal possibility).

In the paper, since we consider the average values of functionals on infinite-dimensional ball, we need a definition of the average value on the ball. Below we give such definitions.

\textbf{ Definition 3}: For $M=\{x|||x||_p\leq R, x(t)\in C[0,1]\}$ where $C[0,1]$ is equipped $p-$norm $||x||_p=(\int_0^1|x(t)|^p\mathrm{d}t)^{\frac{1}{p}}$ for $p\geq1$ and $Y=f(x)=\int_0^1g(x(t))\mathrm{d}t$ where $g(x)$ is a continuous function, and $p=\frac{p_0}{q_0}$ where $p_0$ is even and $(p_0,q_0)=1$,
 we can define the average value of $Y$ on $M$ as
\begin{equation}
EY=\lim_{n\rightarrow \infty}\frac{\int_{M_n} \frac{1}{n}\sum_{k=1}^ng(x_k)\mathrm{d}v_n}{\int_{M_n} \mathrm{d}v_n}
\end{equation}
where $x_k=x(\frac{k}{n})$, $M_n=\{(x_1,\cdots,x_n)|x_1^p+\cdots+x_n^p\leq nR^p\}$ and $\mathrm{d}v_n=\mathrm{d}x_1\cdots \mathrm{d}x_n$ is the volume element of $M_n$. If the limitation exists and is finite or infinite, we call it the average value of functional $f$ on $M$. We often use $EY$ to denote the average value of functional $Y=f(x)$.

\textbf{ Definition 4}: For $M^+=\{x|||x||_p\leq R, x(t)\geq0, x(t)\in C[0,1]\}$ where $C[0,1]$ is equipped $p-$norm and $Y=f(x)=\int_0^1g(x(t))\mathrm{d}t$ where $g(x)$ is a continuous function, and $p\geq 1$ or specially $p=\frac{p_0}{q_0}$ where $p_0$ is odd and $(p_0,q_0)=1$,
 we can define the average value of $Y$ on $M^+$ as
\begin{equation}
EY=\lim_{n\rightarrow \infty}\frac{\int_{M_n^+} \frac{1}{n}\sum_{k=1}^ng(x_k)\mathrm{d}v_n}{\int_{M_n^+} \mathrm{d}v_n}
\end{equation}
where $x_k=x(\frac{k}{n})$, $M_n^+=\{(x_1,\cdots,x_n)|x_1^p+\cdots+x_n^p\leq nR^p, x_k\geq0, k=1,\cdots, n\}$ and $\mathrm{d}v_n=\mathrm{d}x_1\cdots \mathrm{d}x_n$ is the volume element of $M_n^+$. If the limitation exists and is finite or infinite, we call it the average value of functional $f$ on $M^+$.

For $m$ variables function $g(x_1,\cdots,x_m)$, the average value of the functional $f(x)=\int_0^1\cdots \int_0^1g(x(t_1),\cdots,x(t_m))\mathrm{d}t_1\cdots\mathrm{d}t_m$ on infinite-dimensional balls can be defined similarly.

In addition, since there doesn't exist the Lebesque measure on $M$, the probability theory based on the lebesque measure also doesn't exist. Therefore, under the rigorous mathematical meaning, we cannot say what is the probability of the functional $Y$ deviating its average value $EY$. However, in order to measure the deviation of $Y$ from $EY$, we can still define the variance $DY$ to do this because the variance $DY$ is also average value of the functional $(Y-EY)^2$, that is, $DY=E(Y-EY)^2=E(Y^2)-E^2Y$.

For the purpose of discussion on nonlinear exchange formula, we need the definition of average value of $h(Y)$ where $h(Y)$ is a continuous or analytic function of $Y$.

\textbf{Definition 5}. Denote $M$ as the previous $M$ or $M^+$, $Y=f(x)=\int_0^1g(x(t))\mathrm{d}t$ where $g(x)$ is a continuous function, then for a continuous function $h(Y)$,
 we can define the average value of $h(Y)$ on $M$ as
\begin{equation}
Eh(Y)=\lim_{n\rightarrow \infty}\frac{\int_{M_n} h(\frac{1}{n}\sum_{k=1}^ng(x_k))\mathrm{d}v_n}{\int_{M_n} \mathrm{d}v_n}
\end{equation}
where $x_k=x(\frac{k}{n})$, and $M_n$ is also previous $M_n$ or $M^+_n$. If the limitation exists and is finite or infinite, we call it the average value of functional $h(Y)$ on $M$.

This paper is outlined as follows. In section 2,  we drive out some probability densities of coordinates of points in infinite-dimensional balls by two ways of analysis and geometry. In section 3, we give the exact values of some infinite-dimensional integrals. Furthermore, we discuss the
concentration without measure, and obtain the nonlinear exchange formulas for infinite-dimensional integrals. In section 4, we give some further results and definitions.  The last section is short conclusion.

\section{The probability densities of the coordinates of points in infinite-dimensional balls}

We first derive several interesting probability densities from a geometrical way based on the consideration in infinite dimensional space. These results have also independent values.

Consider the continuous functions space $C[0,1]$ and define some norms such as $||x||_0=\max_{t\in[0,1]}|x(t)|$, and $||x||_p=(\int_0^1|x(t)|^p\mathrm{d}t)^{\frac{1}{p}}$ for $p\geq1$. For $p=\frac{p_0}{q_0}$ where $p_0$ is even and $(p_0,q_0)=1$, we consider the whole ball
$M=\{x|||x||_p\leq R\}$, while for $p\geq1$ being a general real number or  specially $p=\frac{p_0}{q_0}$ where  $p_0$ is odd and $(p_0,q_0)=1$, we only consider the "first quadrant" of $M$, that is
$M^+=\{x|x\geq0,||x||_p\leq R\}$.

The following lemma is important.

\textbf{Lemma 1}([33]). The following generalized Dirichlet formula holds
\begin{equation}
\int\cdots\int_{B^+}x_1^{p_1-1}x_2^{p_2-1}\cdots x_n^{p_n-1}\mathrm{d}x_1\cdots\mathrm{d}x_n=
\frac{1}{2^n}\frac{\Gamma(\frac{p_1}{2})\cdots \Gamma(\frac{p_n}{2})}{\Gamma(1+\frac{p_1+\cdots+p_n}{2})},
\end{equation}
where $p_i>0$ for $i=1,\cdots,n$ and $B^+=\{(x_1,\cdots,x_n)|x_1^2+\cdots+x_n^2\leq1, x_k\geq0, k=1,\cdots, n\}$.

Next we give the following  results.

\textbf{Theorem 1}(Version of analysis). For the set $M_n=\{(x_1,\cdots,x_n)|x_1^p+\cdots+x_n^p\leq nR^p\}$ where $p=\frac{p_0}{q_0}$ and $p_0$ is even and $(p_0,q_0)=1$, when we suppose that the points in  $M_n$ are taken by equiprobability, the density of every coordinate $x_k$ of $x$ as a random variable  is given  by
\begin{equation}
\rho_n(x_k)
=\frac{p\Gamma(1+\frac{n}{p})}{2Rn^{\frac{1}{p}}\Gamma(\frac{1}{p})\Gamma(1+\frac{n-1}{p})}
(1-\frac{x_k^p}{nR^p})^{\frac{n-1}{p}}.
\end{equation}
In particular, the limitation of $\rho_n(x)$ as $n$ tending to infinity is given by
\begin{equation}
\rho(x)=\frac{1}{2R\Gamma(\frac{1}{p})p^{\frac{1}{p}-1}}\mathrm{e}^{-\frac{x^p}{pR^p}}, x\in(-\infty, +\infty).
\end{equation}
In general, for any $k$ distinct coordinates $x_{i_1},\cdots, x_{i_k}$, where $({i_1},\cdots, {i_k})\subset\{1,\cdots,n\}$, their union density  is given by
\begin{equation}
\rho_n(x_{i_1},\cdots, x_{i_k})
=\frac{p^k\Gamma(1+\frac{n}{p})}{2^kR^kn^{\frac{k}{p}}\Gamma^k(\frac{1}{p})\Gamma(1+\frac{n-k}{p})}
(1-\frac{x_{i_1}^p+\cdots + x_{i_k}^p}{nR^p})^{\frac{n-k}{p}},
\end{equation}
and the limitation of $\rho_n(x_{i_1},\cdots, x_{i_k})$ as $n$ approaching to  infinity is
\begin{equation}
\rho(x_{i_1},\cdots, x_{i_k})=\frac{1}{2^kR^k\Gamma^k(\frac{1}{p})p^{\frac{k}{p}-k}}\mathrm{e}^{-\frac{x_{i_1}^p+\cdots + x_{i_k}^p}{pR^p}}, x_{i_j}\in(-\infty, +\infty), j=1,\cdots,k,
\end{equation}
that is,
\begin{equation}
\rho(x_{i_1},\cdots, x_{i_k})=\rho(x_{i_1})\cdots\rho( x_{i_k}),
\end{equation}
which means that when $n$ tends to infinity, any finite coordinates of point in $M_n$ as random variables are independent.

\textbf{Proof}. By symmetry, we only consider the density of $(x_1,\cdots,x_k)$. Denote $M'_n=\{(x_{k+1},\cdots,x_n)|x_{k+1}^p+\cdots+x_n^p\leq nR^p-x_1^p-\cdots-x^p_k\}$.  According to the assumption of equiprobability, we have
\begin{equation}
\rho_n(x_1,\cdots,x_k)=\frac{\int_{M'_n}\mathrm{d}x_{k+1}\cdots\mathrm{d}x_n}{\int_{M_n}\mathrm{d}x_1\cdots\mathrm{d}x_n}.
\end{equation}
Further, by taking the transformation $x_k=Rn^{\frac{1}{p}}y_k^{\frac{2}{p}}$, we have from the lemma 1,
\begin{equation}
\rho_n(x_1,\cdots,x_k)
=\frac{p^k\Gamma(1+\frac{n}{p})}{2^kR^k\Gamma^k(\frac{1}{p})n^{\frac{k}{p}}\Gamma(1+\frac{n-k}{p})}
(1-\frac{x_1^p+\cdots+x^p_k}{nR^p})^{\frac{n-k}{p}}.
\end{equation}
Taking the limitation of $n$ approaching to $+\infty$, and using the Stirling's asymptotic formula of Gamma function,  $\Gamma(s)\sim \sqrt{2\pi}s^{s-\frac{1}{2}}\mathrm{e}^{-s}$, we get,
\begin{equation*}
\rho(x_1,\cdots,x_k)=\lim_{n\rightarrow+\infty}\frac{p^k\sqrt{2\pi}\sqrt{\frac{n}{p}}
(\frac{n}{p})^{\frac{n}{p}}\mathrm{e}^{-\frac{n}{p}}}
{2^kR^k\Gamma^k(\frac{1}{p})\sqrt{2\pi}\sqrt{\frac{n-k}{p}}
(\frac{n-k}{p})^{\frac{n-k}{p}}n^{\frac{k}{p}}\mathrm{e}^{-\frac{n-k}{p}}}
(1-\frac{x_1^p+\cdots+x_k^p}{nR^p})^{\frac{n-k}{p}}
\end{equation*}
\begin{equation*}
=\lim_{n\rightarrow+\infty}\frac{p^k\sqrt{\frac{n}{n-k}}
(\frac{n}{n-k})^{\frac{k}{p}}}
{2^kR^k\Gamma^k(\frac{1}{p})
(\frac{n-k}{n})^{\frac{n}{p}}\mathrm{e}^{\frac{k}{p}}}
(1-\frac{x_1^p+\cdots+x_k^p}{nR^p})^{\frac{n-k}{p}}
\end{equation*}
\begin{equation}
=\frac{1}{2^kR^kp^{\frac{k}{p}-k}\Gamma^k(\frac{1}{p})}\mathrm{e}^{-\frac{x_1^p+\cdots+x_k^p}{pR^p}}.
\end{equation}
 The proof is completed.

Now we consider the infinite-dimensional ball $M=\{x|\int_0^1x^p(t)\mathrm{d}t\leq R^p\}$ in $C[0,1]$ with $p-$norm. Although a uniform distribution mathematically does not exist on the ball $M$ because
the dimension on $M$ is infinite, the density of coordinates of points in $M$ does exist! Fortunately, we needn't this uniform distribution to derive our the result. What we only need is a limit procedure from finite dimension to infinite dimension so that we can avoid the trouble of nonexistence of uniform distribution. The following is the version of geometry of theorem 1 under the meaning of definition 2. In other words, this is just a probability "explanation" in formal.

\textbf{Theorem 1}(Version of geometry). For the ball $M=\{x|\int_0^1x^p(t)\mathrm{d}t\leq R^p\}$ in $C[0,1]$ where $p=\frac{p_0}{q_0}$ and $p_0$ is even and $(p_0,q_0)=1$, when we suppose that the points in the ball $M$ are taken by equiprobability under the meaning of definition 2, the density of $x(t)$ as a random variable for fixed $t$ is given  by
\begin{equation}
\rho(x)=\frac{1}{2R\Gamma(\frac{1}{p})p^{\frac{1}{p}-1}}\mathrm{e}^{-\frac{x^p}{pR^p}}, x\in(-\infty, +\infty).
\end{equation}
In general, for any $k$ distinct coordinates $x_{t_1},\cdots, x_{t_k}$, where $0<t_1<\cdots< t_k<1$, their union density  is given by
\begin{equation}
\rho(x_{1},\cdots, x_{k})=\frac{1}{2^kR^k\Gamma^k(\frac{1}{p})p^{\frac{k}{p}-k}}\mathrm{e}^{-\frac{x_{1}^p+\cdots + x_{k}^p}{pR^p}}, x_{j}\in(-\infty, +\infty), j=1,\cdots,k,
\end{equation}
that is,
\begin{equation}
\rho(x_{1},\cdots, x_{k})=\rho(x_{1})\cdots\rho( x_{k}),
\end{equation}
which means that  any finite coordinates of point $X$ in $M$ as random variables are independent.

\textbf{Proof}. Firstly, by discretization(it is reasonable by the continuity of $x(t)$), we have $M_n=\{(x_1,\cdots,x_n)|x_1^p+\cdots+x_n^p\leq nR^p\}$ where
$x_k=x(\frac{k}{n})$. We direct compute the density $\rho_n(x_1,\cdots,x_k)$ of $(x_1,\cdots,x_k)$ as the random variables in the ball $M_n$. According to the assumption of equiprobability and the version of analysis of the theorem 1, the theorem is proven.

This result is different to the finite case essentially. In finite ball, coordinates are not independent each other since they are constrained on the ball and there exists a certain relation. But in infinite dimensional ball, for any finite number of coordinates of point in the ball, from its discretization $M_n=\{(x_1,\cdots,x_n)|x_1^p+\cdots+x_n^p\leq nR^p\}$, we can easily see that as $n$ tending to infinity, the radius tends to also infinity, so, for any finite number of coordinates such as $x_1,\cdots,x_k$, their value ranges will become the whole $n$ dimensional space $R^n$. This means that the constraint has disappeared and hence these coordinates are really independent. In other words, essentially, the ball $M$ contains all finite dimensional linear spaces $R^n$ for any positive integer $n$.

Similarly,  we have the following theorems.

\textbf{Theorem 2}(Version of analysis). For  the "first quadrant" $M^+_n=\{(x_1,\cdots,x_n)|x_1^p+\cdots+x_n^p\leq nR^p, x_k\geq 0, k=1,\cdots,n\}$ where $p$ is a general real number and $p\geq 1$ or specially $p=\frac{p_0}{q_0}$ where $p_0$ is odd and $(p_0,q_0)=1$, when we suppose that the points in $M_n^+$ are taken by equiprobability, the density of every coordinate $x_k$ of $x$ as a random variable  is given by
\begin{equation}
\rho_n(x_k)
=\frac{p\Gamma(1+\frac{n}{p})}{Rn^{\frac{1}{p}}\Gamma(\frac{1}{p})\Gamma(1+\frac{n-1}{p})}
(1-\frac{x_k^p}{nR^p})^{\frac{n-1}{p}}.
\end{equation}
In particular, the limitation of $\rho_n(x)$ as $n$ tending to infinity is given by
\begin{equation}
\rho(x)=\frac{1}{R\Gamma(\frac{1}{p})p^{\frac{1}{p}-1}}\mathrm{e}^{-\frac{x^p}{pR^p}}, x\in[0, +\infty).
\end{equation}
In general, for any $k$ distinct coordinates $x_{i_1},\cdots, x_{i_k}$, where $({i_1},\cdots, {i_k})\subset\{1,\cdots,n\}$, their union density  is given by
\begin{equation}
\rho_n(x_{i_1},\cdots, x_{i_k})
=\frac{p^k\Gamma(1+\frac{n}{p})}{R^kn^{\frac{k}{p}}\Gamma^k(\frac{1}{p})\Gamma(1+\frac{n-k}{p})}
(1-\frac{x_{i_1}^p+\cdots + x_{i_k}^p}{nR^p})^{\frac{n-k}{p}},
\end{equation}
and the limitation of $\rho_n(x_{i_1},\cdots, x_{i_k})$ as $n$ approaching to  infinity is
\begin{equation}
\rho(x_{i_1},\cdots, x_{i_k})=\frac{1}{R^k\Gamma^k(\frac{1}{p})p^{\frac{k}{p}-k}}\mathrm{e}^{-\frac{x_{i_1}^p+\cdots + x_{i_k}^p}{pR^p}}, x_{i_j}\in[0, +\infty), j=1,\cdots,k,
\end{equation}
that is,
\begin{equation}
\rho(x_{i_1},\cdots, x_{i_k})=\rho(x_{i_1})\cdots\rho( x_{i_k}),
\end{equation}
which means that when $n$ tends to infinity, any finite coordinates of point in $M_n$ as random variables are independent.

\textbf{Proof}. According to the assumption of equiprobability, we have from the lemma 1,
\begin{equation*}
\rho_n(x_1,\cdots,x_k)=\frac{\int\cdots\int_{M'^+_n}\mathrm{d}x_{k+1}\cdots\mathrm{d}
x_n}{\int_{M^+_n}\mathrm{d}x_1\cdots\mathrm{d}x_n}
\end{equation*}
\begin{equation}
=\frac{p\Gamma(1+\frac{n}{p})}{\Gamma(\frac{1}{p})\Gamma(1+\frac{n-1}{p})(nR^p-x^p_1)^{\frac{1}{p}}}
(1-\frac{x_1^p}{nR^p})^{\frac{n}{p}},
\end{equation}
where $M'^+_n=\{(x_{k+1},\cdots,x_n)|x_{k+1}^p+\cdots+x_n^p\leq nR^p-x_1^p-\cdots-x_k^p, x_k\geq 0, k=1,\cdots,n\}$.
Taking the limitation of $n$ approaching to $+\infty$, and using the Stirling's asymptotic formula of Gamma function,  we get,
\begin{equation*}
\rho(x_1,\cdots,x_k)=\lim_{n\rightarrow+\infty}\frac{p^k\sqrt{2\pi}\sqrt{\frac{n}{p}}
(\frac{n}{p})^{\frac{n}{p}}\mathrm{e}^{-\frac{n}{p}}}
{2^kR^k\Gamma^k(\frac{1}{p})\sqrt{2\pi}\sqrt{\frac{n-k}{p}}
(\frac{n-k}{p})^{\frac{n-k}{p}}n^{\frac{k}{p}}\mathrm{e}^{-\frac{n-k}{p}}}
(1-\frac{x_1^p+\cdots+x_k^p}{nR^p})^{\frac{n-k}{p}}
\end{equation*}
\begin{equation*}
=\lim_{n\rightarrow+\infty}\frac{p^k\sqrt{\frac{n}{n-k}}
(\frac{n}{n-k})^{\frac{k}{p}}}
{R^k\Gamma^k(\frac{1}{p})
(\frac{n-k}{n})^{\frac{n}{p}}\mathrm{e}^{\frac{k}{p}}}
(1-\frac{x_1^p+\cdots+x_k^p}{nR^p})^{\frac{n-k}{p}}
\end{equation*}
\begin{equation}
=\frac{1}{R^kp^{\frac{k}{p}-k}\Gamma^k(\frac{1}{p})}\mathrm{e}^{-\frac{x_1^p+\cdots+x_k^p}{pR^p}}.
\end{equation}
  The proof is completed.

  Similar to the theorem 1, we have the version of geometry of the theorem 2.

  \textbf{Theorem 2}(Version of geometry). For  the "first quadrant" $M^+=\{x|x(t)\geq0, \int_0^1x^p(t)\mathrm{d}t\leq R^p\}$ of $M$ where $p$ is a general real number and $p\geq 1$ or specially $p=\frac{p_0}{q_0}$ where $p_0$ is odd and $(p_0,q_0)=1$, when we suppose that the points in $M^+$ are taken by equiprobability, the density of $x(t)$ as a random variable on $M^+$ for fixed $t$ is given by
\begin{equation}
\rho(x)=\frac{1}{Rp^{\frac{1}{p}-1}\Gamma(\frac{1}{p})}\mathrm{e}^{-\frac{x^p}{pR^p}}, x\in[0, +\infty).
\end{equation}
In general, for any $k$ distinct coordinates $x_{t_1},\cdots, x_{t_k}$, where $0<t_1<\cdots< t_k<1$, their union density  is given by
\begin{equation}
\rho(x_{1},\cdots, x_{k})=\frac{1}{R^k\Gamma^k(\frac{1}{p})p^{\frac{k}{p}-k}}\mathrm{e}^{-\frac{x_{1}^p+\cdots + x_{k}^p}{pR^p}}, x_{j}\in[0, +\infty), j=1,\cdots,k,
\end{equation}
that is,
\begin{equation}
\rho(x_{1},\cdots, x_{k})=\rho(x_{1})\cdots\rho( x_{k}),
\end{equation}
which means that any finite coordinates of point $X$ in $M$ as random variables are independent.

\textbf{Remark 2}. If we take some special values of $p$ and suitable variable transformation, we will gain some interesting and important probability distributions.  When $p$ is even, the density looks like a normal distribution, and thus we call it high order normal distribution or normal-like distribution.  If $R=1$, a simple form is
\begin{equation}
\rho(x)=\frac{1}{2\Gamma(\frac{1}{p})p^{\frac{1}{p}-1}}\mathrm{e}^{-\frac{x^p}{p}}, x\in(-\infty, +\infty).
\end{equation}
Further, for example, if $p=2$, we get the standard normal distribution
\begin{equation}
\rho(x)=\frac{1}{\sqrt{2\pi}}\mathrm{e}^{-\frac{x^2}{2}}, x\in(-\infty, +\infty),
\end{equation}
which gives the G\^{a}teaux and L\'{e}vy's result [4].
If $p=4$, we get a  4-order normal distribution
\begin{equation}
\rho(x)=\frac{\sqrt2}{\Gamma(\frac{1}{4})}\mathrm{e}^{-\frac{x^4}{4}}, x\in(-\infty, +\infty).
\end{equation}

 When $p$ is odd, the density looks like an exponent distribution, and thus we call it high order exponent distribution or exponent-like distribution.  For example, if $p=1$ and $R=\lambda^{-\frac{1}{p}} $, we get the usual exponent distribution
\begin{equation}
\rho(x)=\frac{1}{\lambda}\mathrm{e}^{-\lambda x}, x\in[0, +\infty).
\end{equation}
If $p=3$ and $R=(3\lambda)^{-\frac{1}{p}} $, we get the 3-order exponent distribution
\begin{equation}
\rho(x)=\frac{3\lambda^{\frac{1}{3}}}{\Gamma(\frac{1}{3})}\mathrm{e}^{-\lambda x^3}, x\in[0, +\infty).
\end{equation}

 By a simple transformation, we can obtain the famous Gamma distribution in statistics. Indeed, we take a transformation
\begin{equation}
Z=\frac{x^p}{p\beta},
\end{equation}
then the density of $Z$ is just
\begin{equation}
\rho(z)=\frac{\beta^{\frac{1}{p}}}{\Gamma(\frac{1}{p})}z^{\frac{1}{p}-1}\mathrm{e}^{-\beta y}.
\end{equation}
Further, taking $\alpha=\frac{1}{p}$ gives the Gamma distribution
\begin{equation}
\rho(z,\alpha,\beta)=\frac{\beta^{\alpha}}{\Gamma(\alpha)}z^{\alpha-1}\mathrm{e}^{-\beta y}.
\end{equation}
 This is a geometrical origin of the Gamma distribution. We can see that this is a rather natural way to derive the Gamma distribution.

 Based on the maximum non-symmetrical entropy principle, we can also derive out some distributions [34,35]. But the above geometrical origin is more natural.

\section{The average values of some functionals and the concentration without measure}

In the section, according to the above results, we study a kind of infinite-dimensional functionals with the form of integral. Our main results are summarized in theorems 3 and 4. The considered infinite-dimensional integrals arise from  the infinite-dimensional probability theory[36]. An elementary and rough introduction on such topic can be seen in [36] in which we give more examples and another way to compute the exact average values of some functionals.

\textbf{Lemma 2}. If $f(x)$ satisfies one of the following two conditions,

$H_1$. (Differential condition):  $f(x)\in C^k[0,+\infty)$, $\lim_{x\rightarrow+\infty}f^{(j)}(x)\mathrm{e}^{-x}=0$ for $j=0,\cdots, k-1$, and there is a positive constant number $A$ such that  $|f^{(k)}(x)|\leq A$;

$H_2$. (Integral condition):
\begin{equation}
\int_0^{+\infty}|f(x)|\mathrm{e}^{-\frac{x}{p}}\mathrm{d}x<+\infty,\int_0^{+\infty}x^2|f(x)|\mathrm{e}^{-\frac{x}{p}}\mathrm{d}x<+\infty;
\end{equation}
then we have
\begin{equation}
\lim_{n\rightarrow +\infty}\int_0^n f(x)(1-\frac{x}{n})^{\frac{n}{p}}\mathrm{d}x=\int_0^{+\infty}f(x)\mathrm{e}^{-\frac{x}{p}}\mathrm{d}x.
\end{equation}
In general, for any finite integer $n_0$, we also have
\begin{equation}
\lim_{n\rightarrow +\infty}\int_0^n f(x)(1-\frac{x}{n})^{\frac{n-n_0}{p}}\mathrm{d}x=\int_0^{+\infty}f(x)\mathrm{e}^{-\frac{x}{p}}\mathrm{d}x.
\end{equation}

\textbf{Proof}. \textbf{Case (i)}. We  prove the general formula. Under the condition $H_1$. From integration by parts, we have
\begin{equation*}
\int_0^n f(x)(1-\frac{x}{n})^{\frac{n-n_0}{p}}\mathrm{d}x=\frac{pn}{n-n_0+p}f(0)+\frac{pn}{n-n_0+p}\int_0^n f'(x)(1-\frac{x}{n})^{\frac{n-n_0}{p}+1}\mathrm{d}x
\end{equation*}
\begin{equation*}
=\frac{pn}{n-n_0+p}f(0)
+\frac{p^2n^2}{(n-n_0+p)(n-n_0+2p)}f'(0)+\cdots+\frac{p^kn^k}{(n-n_0+p)\cdots(n-n_0+kp)}f^{(k-1)}(0)
\end{equation*}
\begin{equation}
+\frac{p^kn^k}{(n-n_0+p)\cdots(n-n_0+kp)}\int_0^n f^{(k)}(x)(1-\frac{x}{n})^{\frac{n-n_0}{p}+k}\mathrm{d}x,
\end{equation}
and
\begin{equation}
\int_0^{+\infty}f(x)\mathrm{e}^{-\frac{x}{p}}\mathrm{d}x=pf(0)+p^2f'(0)+\cdots+p^kf^{(k-1)}(0)
+p^k\int_0^{+\infty}f^{(k)}(x)\mathrm{e}^{-\frac{x}{p}}\mathrm{d}x.
\end{equation}
It is easy to see that we only need to prove
\begin{equation}
\lim_{n\rightarrow +\infty}\int_0^n f^{(k)}(x)(1-\frac{x}{n})^{\frac{n-n_0}{p}+k}\mathrm{d}x=\int_0^{+\infty}f^{(k)}(x)\mathrm{e}^{-\frac{x}{p}}\mathrm{d}x.
\end{equation}
Indeed, by mean value theorem of integral, we have
\begin{equation*}
\lim_{n\rightarrow +\infty}|\int_0^n f^{(k)}(x)((1-\frac{x}{n})^{\frac{n-n_0}{p}+k}-{e}^{-\frac{x}{p}})\mathrm{d}x|
\end{equation*}
\begin{equation*}
=\lim_{n\rightarrow +\infty}|f^{(k)}(\xi_n)||\int_0^n ((1-\frac{x}{n})^{\frac{n-n_0}{p}+k}-{e}^{-\frac{x}{p}})\mathrm{d}x|
\end{equation*}
\begin{equation*}
\leq A\lim_{n\rightarrow +\infty}|\int_0^n \{(1-\frac{x}{n})^{\frac{n-n_0}{p}+k}-{e}^{-\frac{x}{p}}\}\mathrm{d}x|
\end{equation*}
\begin{equation}
= A\lim_{n\rightarrow +\infty}|\frac{np}{n-n_0+(k+1)p}-p+p\mathrm{e}^{-\frac{n}{p}}|=0,
\end{equation}
where $\xi_n\in(0,n)$. The lemma is proven under the first condition (i).

\textbf{Case (ii)}.  Under the condition $H_2$. For $0\leq x\leq n$, we know
\begin{equation}
1+\frac{x}{n}\leq \mathrm{e}^{\frac{x}{n}}\leq\frac{1}{1-\frac{x}{n}},
\end{equation}
and then for $p>0$
\begin{equation}
(1+\frac{x}{n})^{\frac{n}{p}}\leq \mathrm{e}^{\frac{x}{p}}, (1-\frac{x}{n})^{\frac{n}{p}}\leq \mathrm{e}^{-\frac{x}{p}}.
\end{equation}
Therefore,
\begin{equation*}
0\leq \mathrm{e}^{-\frac{x}{p}}-(1-\frac{x}{n})^{\frac{n}{p}}=\mathrm{e}^
{-\frac{x}{p}}\{1-\mathrm{e}^{\frac{x}{p}}(1-\frac{x}{n})^{\frac{n}{p}}\}
\end{equation*}
\begin{equation*}
\leq \mathrm{e}^
{-\frac{x}{p}}\{1-(1-\frac{x^2}{n^2})^{\frac{n}{p}}\}
\end{equation*}
\begin{equation*}
=\mathrm{e}^
{-\frac{x}{p}}(1-(1-\frac{x^2}{n^2})^{\frac{1}{p}})\{1+(1-\frac{x^2}{n^2})^
{\frac{1}{p}}+\cdots+(1-\frac{x^2}{n^2})^{\frac{n-1}{p}}\}
\end{equation*}
\begin{equation}
\leq n\mathrm{e}^
{-\frac{x}{p}}(1-(1-\frac{x^2}{n^2})^{\frac{1}{p}}).
\end{equation}
In addition, we have
\begin{equation}
\lim_{n\rightarrow +\infty}\frac{1-(1-\frac{x^2}{n^2})^{\frac{1}{p}}}{\frac{x^2}{pn^2}}=1,
\end{equation}
and hence for arbitrary $0<\epsilon<1$, there exists constant $N_1$, such that for $n>N_1$,
\begin{equation}
|\frac{1-(1-\frac{x^2}{n^2})^{\frac{1}{p}}}{\frac{x^2}{pn^2}}-1|<\epsilon.
\end{equation}
by which, we have for $n>N_1$,
\begin{equation*}
|\int_0^n f(x)({e}^{-\frac{x}{p}}-(1-\frac{x}{n})^{\frac{n}{p}})\mathrm{d}x|
\end{equation*}
\begin{equation*}
\leq\frac{1}{pn}\int_0^n x^2|f(x)|\mathrm{e}^
{-\frac{x}{p}}\mathrm{d}x+\frac{1}{pn}\int_0^n x^2|f(x)|\mathrm{e}^
{-\frac{x}{p}}|\frac{1-(1-\frac{x}{n})^{\frac{n}{p}}}{\frac{x^2}{pn^2}}-1|\mathrm{d}x
\end{equation*}
\begin{equation}
\leq\frac{1}{pn}\int_0^{+\infty} x^2|f(x)|\mathrm{e}^
{-\frac{x}{p}}\mathrm{d}x+\frac{\epsilon}{pn}\int_0^{+\infty} x^2|f(x)|\mathrm{e}^
{-\frac{x}{p}}\mathrm{d}x.
\end{equation}
Since there exists an enough large number $N_2$ such that for $n>N_2$, we have
\begin{equation}
\frac{1}{pn}\int_0^{+\infty} x^2|f(x)|\mathrm{e}^
{-\frac{x}{p}}\mathrm{d}x<\frac{\epsilon}{2},
\end{equation}
 taking $N=\max\{N_1,N_2\}$, for $n>N$, we get
 \begin{equation}
|\int_0^n f(x)({e}^{-\frac{x}{p}}-(1-\frac{x}{n})^{\frac{n}{p}})\mathrm{d}x|<\frac{\epsilon}{2}+\frac{\epsilon^2}{2}<\epsilon.
\end{equation}
Since $\lim_{n\rightarrow +\infty}(1-\frac{x}{n})^{-\frac{n_0}{p}}=1$ and $\int_0^{+\infty}|f(x)|\mathrm{e}^{-\frac{x}{p}}\mathrm{d}x<+\infty$, there exists an enough number $N_3$ such that for $n>N_3$,
\begin{equation}
|\int_0^n f(x)(1-\frac{x}{n})^{\frac{n}{p}}\{(1-\frac{x}{n})^{-\frac{n_0}{p}}-1\}\mathrm{d}x|
\leq\int_0^n |f(x)|\mathrm{e}^{-\frac{x}{p}}|(1-\frac{x}{n})^{-\frac{n_0}{p}}-1|\mathrm{d}x|<\epsilon.
\end{equation}
Thus, we have for $n>\max\{N_1,N_2,N_3\}$
\begin{equation*}
|\int_0^n f(x)({e}^{-\frac{x}{p}}-(1-\frac{x}{n})^{\frac{n-n_0}{p}})\mathrm{d}x|
\end{equation*}
 \begin{equation*}
\leq|\int_0^n f(x)({e}^{-\frac{x}{p}}-(1-\frac{x}{n})^{\frac{n}{p}})\mathrm{d}x|
\end{equation*}
\begin{equation}
+|\int_0^n f(x)(1-\frac{x}{n})^{\frac{n}{p}}\{(1-\frac{x}{n})^{-\frac{n_0}{p}}-1\}\mathrm{d}x|
<2\epsilon.
\end{equation}
The lemma is proven.

It is easy to generalize the lemma 2 to the case of $f(x)$ depending on $n$.

\textbf{Lemma 3}. Suppose  $f_n(x)$  satisfies one of the following two conditions:

$H_1$. (Differential condition):  $\lim_{n\rightarrow+\infty}f^{(j)}_n(x)=f^{(j)}(x)$ for $j=0,\cdots,k$, $f_n(x)\in C^k[0,+\infty)$, $\lim_{x\rightarrow+\infty}f_n^{(j)}(x)\mathrm{e}^{-x}=0$ for $j=0,\cdots, k-1$, and there is a positive constant number $A$ such that  $|f_n^{(k)}(x)|\leq A$;

$H_2$. (Integral condition): Suppose that  $\lim_{n\rightarrow+\infty}f_n(x)=f(x)$ uniformly for $x$, and
\begin{equation}
\int_0^{+\infty}|f(x)|\mathrm{e}^{-\frac{x}{p}}\mathrm{d}x<+\infty,
\int_0^{+\infty}x^2|f(x)|\mathrm{e}^{-\frac{x}{p}}\mathrm{d}x<+\infty.
\end{equation}
Then we have for any finite integer $n_0$,
\begin{equation}
\lim_{n\rightarrow +\infty}\int_0^n f_n(x)(1-\frac{x}{n})^{\frac{n-n_0}{p}}\mathrm{d}x=\int_0^{+\infty}f(x)\mathrm{e}^{-\frac{x}{p}}\mathrm{d}x.
\end{equation}

\textbf{Proof}. The proof in the condition $H_1$ is similar completely to the lemma 2, only replacing $f(x)$ by $f_n(x)$.
In the condition $H_2$, we have
\begin{equation*}
|\int_0^n \{f(x){e}^{-\frac{x}{p}}-f_n(x)(1-\frac{x}{n})^{\frac{n-n_0}{p}}\}\mathrm{d}x|
\end{equation*}
\begin{equation}
\leq|\int_0^n f(x)({e}^{-\frac{x}{p}}-(1-\frac{x}{n})^{\frac{n-n_0}{p}})\mathrm{d}x|+|\int_0^n (f(x)-f_n(x))(1-\frac{x}{n})^{\frac{n-n_0}{p}}\mathrm{d}x|.
\end{equation}
The first term is the same as that in lemma 2, and the second term tends to zero by the assumption of uniform convergence of $f_n(x)$ and $\int_0^n(1-\frac{x}{n})^{\frac{n-n_0}{p}}\mathrm{d}x=\frac{p}{n-n_0+p}$. The proof is complete.

Now, for convenience, we define a symmetrization operator $S$ for $m$ variables function $g(x_1,\cdots,x_m)$ by
\begin{equation}
(Sg)(x_1,\cdots,x_m)=\sum_{x\in S(m)}g(x),
\end{equation}
where
\begin{equation}
S(m)=\{x_1,-x_1\}\times\{x_2,-x_2\}\times\cdots\times\{x_m,-x_m\}.
\end{equation}
For example,
\begin{equation}
S(1)=\{x_1,-x_1\}, (Sg)(x)=g(x)+g(-x)
\end{equation}
\begin{equation}
S(2)=\{(x_1,x_2), (x_1,-x_2),(-x_1,x_2), (-x_1,-x_2)\},
\end{equation}
\begin{equation}
(Sg)(x_1,x_2)=g(x_1,x_2)+g(x_1,-x_2)+g(-x_1,x_2)+g(-x_1,-x_2)\}.
\end{equation}

\textbf{Theorem 3}. Denote for $m$ variables function $g$
\begin{equation*}
G(z_1,\cdots,z_{m-1},z)=(Sg)(Rz^{\frac{1}{p}}_1,\cdots,R(z-z_1-\cdots-z_{m-1})^{\frac{1}{p}})
\end{equation*}
\begin{equation}
\times\{z_1(z_2-z_1)\cdots(z-z_1-\cdots-z_{m-1})\}^{\frac{1}{p}-1},
\end{equation}
and
\begin{equation}
f_n(z)=\int_0^n\cdots\int_0^nG(z_1,\cdots,z_{m-1},z))\mathrm{d}z_1\cdots\mathrm{d}z_{m-1},
\end{equation}
\begin{equation}
f(z)=\int_0^{+\infty}\cdots\int_0^{+\infty}G(z_1,\cdots,z_{m-1},z))\mathrm{d}z_1\cdots\mathrm{d}z_{m-1}.
\end{equation}
Suppose that $g$ satisfies one of the following two conditions

$(H_1)$. (Differential condition):  $\lim_{n\rightarrow+\infty}f^{(j)}_n(x)=f^{(j)}(x)$ for $j=0,\cdots,k$, $f_n(x)\in C^k[0,+\infty)$, $\lim_{x\rightarrow+\infty}f_n^{(j)}(x)\mathrm{e}^{-x}=0$ for $j=0,\cdots, k-1$, and there is a positive constant number $A$ such that  $|f_n^{(k)}(x)|\leq A$;

$(H_2)$. (Integral condition): Suppose that  $\lim_{n\rightarrow+\infty}f_n(x)=f(x)$ uniformly for $x$, and
\begin{equation}
\int_0^{+\infty}|f(x)|\mathrm{e}^{-\frac{x}{p}}\mathrm{d}x<+\infty,
\int_0^{+\infty}x^2|f(x)|\mathrm{e}^{-\frac{x}{p}}\mathrm{d}x<+\infty.
\end{equation}
Take the ball $M=\{x|\int_0^1x^p(t)\mathrm{d}t\leq R^p\}$ in $C[0,1]$ with norm $||x||_p$ when  $p=\frac{p_0}{q_0}$ where $p_0$ is even and $(p_0,q_0)=1$. Then for functional
\begin{equation}
Y=\int_0^1\cdots \int_0^1g(x(t_1),\cdots, x(t_m))\mathrm{d}t_1\cdots\mathrm{d}t_m,
\end{equation}
the average value of $Y$ on $M$ satisfies
\begin{equation}
EY=\frac{1}{2^m\Gamma^m(\frac{1}{p})p^{\frac{m}{p}-m}}\int_{-\infty}^{+\infty}\cdots\int_{-\infty}^{+\infty}
g(Rx_1,\cdots,Rx_m)\mathrm{e}^{-\frac{x_1^p}{p}-\cdots-\frac{x_m^p}{p}}\mathrm{d}x_1\cdots\mathrm{d}x_m.
\end{equation}

\textbf{Proof}. By discretization of $Y$, we have
\begin{equation}
Y_n=\frac{1}{n^m}\sum_{i_1=1}^{n}...\sum_{i_m=1}^{n}g(x_{i_1},\cdots,x_{i_m}).
\end{equation}
By symmetry and consideration in combinatorics,  we get
\begin{equation}
\lim_{n\rightarrow+\infty}EY_n=Eg(x_1,\cdots,x_m).
\end{equation}
On the other hand, we know that
\begin{equation}
Eg(x_1,\cdots,x_m)=\lim_{n\rightarrow \infty}\frac{\int_{M_n}g(x_1,\cdots,x_m)\mathrm{d}x_{1}\cdots\mathrm{d}x_n}{\int_{M_n}\mathrm{d}x_1\cdots\mathrm{d}x_n}.
\end{equation}
According to the theorem 1, we have
\begin{equation}
Eg(x_1,\cdots,x_m)=\lim_{n\rightarrow \infty}\int_{M'} g(x_1,\cdots,x_m)\rho_n(x_1,\cdots,x_m)\mathrm{d}v_m,
\end{equation}
where $M'=\{(x_1,\cdots,x_m)|x_1^p+\cdots+x_m^p\leq nR^p\}$. Therefore, we need to prove
\begin{equation*}
\lim_{n\rightarrow \infty}\int_{M'} g(x_1,\cdots,x_m)\frac{p^k\Gamma(1+\frac{n}{p})}{2^mR^mn^{\frac{m}{p}}\Gamma^m(\frac{1}{p})\Gamma(1+\frac{n-m}{p})}
(1-\frac{x_{1}^p+\cdots + x_{m}^p}{nR^p})^{\frac{n-m}{p}}\mathrm{d}v_m,
\end{equation*}
\begin{equation}
=\frac{1}{2^mR^m\Gamma^m(\frac{1}{p})p^{\frac{m}{p}-m}}\int_{-\infty}^{+\infty}\cdots\int_{-\infty}^{+\infty}
g(x_1,\cdots,x_m)\mathrm{e}^{-\frac{x_{1}^p+\cdots + x_{m}^p}{nR^p}}\mathrm{d}x_1\cdots\mathrm{d}x_m.
\end{equation}
It is equivalent to prove
\begin{equation*}
\lim_{n\rightarrow \infty}\int_{M'} g(x_1,\cdots,x_m)\frac{p^m\sqrt{\frac{n}{n-m}}
(\frac{n}{n-m})^{\frac{m}{p}}}
{2^mR^m\Gamma^m(\frac{1}{p})
(\frac{n-m}{n})^{\frac{n}{p}}\mathrm{e}^{\frac{m}{p}}}
(1-\frac{x_1^p+\cdots+x_m^p}{nR^p})^{\frac{n-m}{p}}\mathrm{d}v_m,
\end{equation*}
\begin{equation}
=\frac{1}{2^mR^m\Gamma^m(\frac{1}{p})p^{\frac{m}{p}-m}}\int_{-\infty}^{+\infty}\cdots\int_{-\infty}^{+\infty}
g(x_1,\cdots,x_m)\mathrm{e}^{-\frac{x_{1}^p+\cdots + x_{m}^p}{nR^p}}\mathrm{d}x_1\cdots\mathrm{d}x_m.
\end{equation}
It is enough to prove
\begin{equation*}
\lim_{n\rightarrow \infty}\int_{M'} g(x_1,\cdots,x_m)
(1-\frac{x_1^p+\cdots+x_m^p}{nR^p})^{\frac{n-m}{p}}\mathrm{d}v_m,
\end{equation*}
\begin{equation}
=\int_{-\infty}^{+\infty}\cdots\int_{-\infty}^{+\infty}
g(x_1,\cdots,x_m)\mathrm{e}^{-\frac{x_{1}^p+\cdots + x_{m}^p}{nR^p}}\mathrm{d}x_1\cdots\mathrm{d}x_m.
\end{equation}
 For the purpose of simplicity, we reduce the integration on the whole ball $M'$ to the first quadrant $M'^+=\{(x_1,\cdots,x_m)|x_1^p+\cdots+x_m^p\leq nR^p, x_j\geq0, 1\leq j\leq m\}$. By the symmetrization operator $S$, we have
\begin{equation*}
\int_{M'} g(x_1,\cdots,x_m)
(1-\frac{x_1^p+\cdots+x_m^p}{nR^p})^{\frac{n-m}{p}}\mathrm{d}v_m
\end{equation*}
\begin{equation}
=\int_{M'^+} (Sg)(x_1,\cdots,x_m)
(1-\frac{x_1^p+\cdots+x_m^p}{nR^p})^{\frac{n-m}{p}}\mathrm{d}v_m.
\end{equation}
Further, we take the transformation $x_j=Ry_j^{\frac{1}{p}}$ for $j=1,\cdots,m$, and hence
\begin{equation*}
\int_{M'^+} (Sg)(x_1,\cdots,x_m)
(1-\frac{x_1^p+\cdots+x_m^p}{nR^p})^{\frac{n-m}{p}}\mathrm{d}v_m.
\end{equation*}
\begin{equation}
\frac{R^m}{p^m}\int_{\Delta_m} (Sg)(Ry^{\frac{1}{p}}_1,\cdots,y^{\frac{1}{p}}_m)(y_1\cdots y_m)^{\frac{1}{p}-1}
(1-\frac{y_1+\cdots+y_m}{n})^{\frac{n-m}{p}}\mathrm{d}y_1\cdots\mathrm{d}y_m,
\end{equation}
where $\Delta_m=\{(y_1,\cdots,y_m)|y_1+\cdots+y_m\leq n, y_j\geq0, 1\leq j\leq m\}$.

By the same transformation and symmetrization, we reduce the integration from $R^n$ to $R^{+n}$
\begin{equation*}
\int_{-\infty}^{+\infty}\cdots\int_{-\infty}^{+\infty}
g(x_1,\cdots,x_m)\mathrm{e}^{-\frac{x_{1}^p+\cdots + x_{m}^p}{nR^p}}\mathrm{d}x_1\cdots\mathrm{d}x_m
\end{equation*}
\begin{equation}
=\frac{R^m}{p^m}\int_{0}^{+\infty}\cdots\int_{0}^{+\infty}
(Sg)(y_1,\cdots,y_m)(y_1\cdots y_m)^{\frac{1}{p}-1}\mathrm{e}^{-\frac{y_{1}+\cdots + y_{m}}{n}}\mathrm{d}y_1\cdots\mathrm{d}y_m.
\end{equation}

Further, by variables transformations $z_j=y_1+\cdots+y_k, 1\leq j\leq m$ and denoting $z=z_m$, we have
\begin{equation*}
\int_{\Delta_m} (Sg)(Ry^{\frac{1}{p}}_1,\cdots,y^{\frac{1}{p}}_m)(y_1\cdots y_m)^{\frac{1}{p}-1}
(1-\frac{y_1+\cdots+y_m}{n})^{\frac{n-m}{p}}\mathrm{d}y_1\cdots\mathrm{d}y_m,
\end{equation*}
\begin{equation*}
=\int_{0}^{n}\cdots\int_{0}^{n}(Sg)(Rz^{\frac{1}{p}}_1,\cdots,R(z-z_1-\cdots-z_{m-1})^{\frac{1}{p}})
\end{equation*}
\begin{equation}
\{z_1(z_2-z_1)\cdots(z-z_1-\cdots-z_{m-1})\}^{\frac{1}{p}-1}
(1-\frac{z}{n})^{\frac{n-m}{p}}\mathrm{d}z_1\cdots\mathrm{d}z_{m-1}\mathrm{d}z,
\end{equation}
and
\begin{equation*}
\int_{0}^{+\infty}\cdots\int_{0}^{+\infty}
(Sg)(y_1,\cdots,y_m)(y_1\cdots y_m)^{\frac{1}{p}-1}\mathrm{e}^{-\frac{y_{1}+\cdots + y_{m}}{n}}\mathrm{d}y_1\cdots\mathrm{d}y_m.
\end{equation*}
\begin{equation*}
=\int_{0}^{+\infty}\cdots\int_{0}^{+\infty}(Sg)(Rz^{\frac{1}{p}}_1,\cdots,R(z-z_1-\cdots-z_{m-1})^{\frac{1}{p}})
\end{equation*}
\begin{equation}
\times\{z_1(z_2-z_1)\cdots(z-z_1-\cdots-z_{m-1})\}^{\frac{1}{p}-1}
\mathrm{e}^{-\frac{z}{p}}\mathrm{d}z_1\cdots\mathrm{d}z_{m-1}\mathrm{d}z,
\end{equation}
Denote $f_n(z)$ and $f(z)$ by
\begin{equation*}
f_n(z)=\int_{0}^{n}\cdots\int_{0}^{n}(Sg)(Rz^{\frac{1}{p}}_1,\cdots,R(z-z_1-\cdots-z_{m-1})^{\frac{1}{p}})
\end{equation*}
\begin{equation}
\times\{z_1(z_2-z_1)\cdots(z-z_1-\cdots-z_{m-1})\}^{\frac{1}{p}-1}
\mathrm{d}z_1\cdots\mathrm{d}z_{m-1},
\end{equation}
\begin{equation*}
f(z)=\int_{0}^{+\infty}\cdots\int_{0}^{+\infty}(Sg)(Rz^{\frac{1}{p}}_1,\cdots,R(z-z_1-\cdots-z_{m-1})^{\frac{1}{p}})
\end{equation*}
\begin{equation}
\times\{z_1(z_2-z_1)\cdots(z-z_1-\cdots-z_{m-1})\}^{\frac{1}{p}-1}
\mathrm{d}z_1\cdots\mathrm{d}z_{m-1}.
\end{equation}
Under the conditions and by the lemma 3, we get the conclusion. The proof is completed.

\textbf{Theorem 4}.  Suppose that $g$ satisfies one of the following two conditions:

$H_1$: Let $g(x)\in C^k(-\infty,+\infty)$. Denote  $r(x)=x^{\frac{1}{p}-1}\{g(Rx^{\frac{1}{p}})+g(-Rx^{\frac{1}{p}})\}$, and  $r(x)\in C^k[0,+\infty)$,  and $\lim_{x\rightarrow+\infty}r^{(j)}(x)\mathrm{e}^{-x}=0$ for $j=0,\cdots, k-1$, and there is a positive number $A$ such that  $|r^{(k)}(x)|\leq A$.

$H_2$: $
\int_0^{+\infty}x^{\frac{1}{p}+1}|g(Rx^{\frac{1}{p}})+g(-Rx^{\frac{1}{p}})|\mathrm{e}^{-\frac{x}{p}}\mathrm{d}x<+\infty$.

Take the ball $M=\{x|\int_0^1x^p(t)\mathrm{d}t\leq R^p\}$ in $C[0,1]$ with norm $||x||_p$ when  $p=\frac{p_0}{q_0}$ where $p_0$ is even and $(p_0,q_0)=1$. Then we have

(i). The average value and variance on $M$ of the functional
\begin{equation}
Y=\int_0^1g(x(t))\mathrm{d}t,
\end{equation}
satisfy
\begin{equation}
EY=\frac{1}{2\Gamma(\frac{1}{p})p^{\frac{1}{p}-1}}\int_{-\infty}^{+\infty}
g(Rx)\mathrm{e}^{-\frac{x^p}{p}}\mathrm{d}x,
\end{equation}
\begin{equation}
DY=0.
\end{equation}

(ii) For the functionals $Y_1,\cdots, Y_n$ with the form
 \begin{equation}
Y_k=\int_{I_k}g(x(t))\mathrm{d}t,
\end{equation}
where $I_{k}$ for $ k=1,\cdots,n$ are subintervals of the interval $[0,1]$, then
\begin{equation}
EY_k=\frac{l(I_k)}{2\Gamma(\frac{1}{p})p^{\frac{1}{p}-1}}\int_{-\infty}^{+\infty}
g(Rx)\mathrm{e}^{-\frac{x^p}{p}}\mathrm{d}x,
\end{equation}
\begin{equation}
DY_k=0,
\end{equation}
where $l(I_k)$ is the length of interval $I_k$.

 (iii) Given a function $h\in C^k(-\infty, +\infty)$, and there is a positive constant $A$ such that $|h^{(k)}(y)|<A$,  we have the nonlinear exchange formula for the  average value of $h(Y)$ on $M$,
\begin{equation}
Eh(Y)=h(EY).
\end{equation}

\textbf{Proof}. By the theorem 3, we can derive the first conclusion. Here we give a proof for the first result by lemma 2. We only consider the first case of $H_1$, while the second condition can be similarly proven according to the case (ii) in lemma 2.

(i) We direct use the theorem 1 and the lemma 2 to give a simple proof.  In fact, by definition 3,  we take the same partition of $t$ for both $M$ and $Y$, we have
\begin{equation}
EY=\lim_{n\rightarrow \infty}\frac{\int_{M_n} \frac{1}{n}\sum_{k=1}^ng(x_k)\mathrm{d}v_n}{\int_{M_n} \mathrm{d}v_n},
\end{equation}
where $M_n$ is the same that in theorem 1. By symmetry and theorem 1,  we know that
\begin{equation*}
EY=\lim_{n\rightarrow \infty}\int_{-Rn^{\frac{1}{p}}}^{Rn^{\frac{1}{p}}}g(x)\rho_n(x)\mathrm{d}x
=\lim_{n\rightarrow \infty}\int_{0}^{Rn^{\frac{1}{p}}}\{g(x)+g(-x)\}\rho_n(x)\mathrm{d}x
\end{equation*}
\begin{equation*}
=\lim_{n\rightarrow \infty}\int_{0}^{Rn^{\frac{1}{p}}}\{g(x)+g(-x)\}\frac{p\Gamma(1+\frac{n}{p})}{2\Gamma(\frac{1}{p})
\Gamma(1+\frac{n-1}{p})(nR^p-x^p)^{\frac{1}{p}}}
(1-\frac{x^p}{nR^p})^{\frac{n}{p}}\mathrm{d}x
\end{equation*}
\begin{equation}
=\lim_{n\rightarrow \infty}\int_{0}^{n}\{g(Rx^{\frac{1}{p}})+g(-Rx^{\frac{1}{p}})\}\frac{x^{\frac{1}{p}-1}\Gamma(1+\frac{n}{p})}{2\Gamma(\frac{1}{p})
\Gamma(1+\frac{n-1}{p})n^{\frac{1}{p}}}
(1-\frac{x}{n})^{\frac{n-1}{p}}\mathrm{d}x.
\end{equation}
Further, by Stirling formula of Gamma function and the lemma 2 in which we take $n_0=1$, we obtain
\begin{equation*}
EY=\frac{1}{2p^{\frac{1}{p}}\Gamma(\frac{1}{p})}\int_{0}^{+\infty}x^{\frac{1}{p}-1}\{g(Rx^{\frac{1}{p}})+g(-Rx^{\frac{1}{p}})\}
\mathrm{e}^{-\frac{x}{p}}\mathrm{d}x
\end{equation*}
\begin{equation*}
=\frac{1}{2p^{\frac{1}{p}-1}\Gamma(\frac{1}{p})}\int_{0}^{+\infty}\{g(Rx)+g(-Rx)\}
\mathrm{e}^{-\frac{x^p}{p}}\mathrm{d}x
\end{equation*}
\begin{equation}
=\frac{1}{2p^{\frac{1}{p}-1}\Gamma(\frac{1}{p})}\int_{-\infty}^{+\infty}g(Rx)
\mathrm{e}^{-\frac{x^p}{p}}\mathrm{d}x.
\end{equation}
In the last step, we use the property of $x^p$ being even function by the value of $p$. Then we have proven the first conclusion. This also means the linearity of the expectation $E$,
\begin{equation}
EY=\int_0^1E(g(x(t))\mathrm{d}t=E(g(x(t))
=\frac{1}{2R\Gamma(\frac{1}{p})p^{\frac{1}{p}-1}}\int_{-\infty}^{+\infty}g(x)\mathrm{e}^{-\frac{x^p}{pR^p}}\mathrm{d}x.
\end{equation}
Similarly,
\begin{equation*}
E(Y^2)=\lim_{n\rightarrow+\infty}E(Y^2_n)=\lim_{n\rightarrow+\infty}\frac{1}{n^2}\sum_{i=1}^n\sum_{j=1}^nE(g(x_i)g(x_j))
\end{equation*}
\begin{equation*}
=\lim_{n\rightarrow+\infty}\{\frac{n-1}{n}E(g(x_1)g(x_2))+\frac{1}{n}E(g^2(x_1))\}
\end{equation*}
\begin{equation}
=E(g(x_1)g(x_2)).
\end{equation}
By the theorem 3 in which we take $m=2$ and $g(x_1,x_2)=g(x_1)g(x_2)$, we get .
\begin{equation}
=E(g(x_1)g(x_2))=E(g(x_1))E(g(x_2))=E^2Y,
\end{equation}
which means  $DY=0$.

\textbf{Note}: Formally, this conclusion is obvious. Indeed, by the two-dimensional measure of the set $\{(t,s)|t=s\}$ as a subset of $[0,1]^2$ being zero, and the independence of $x(t)$ and $x(s)$ for $t\neq s$, we have from the theorem 1,
\begin{equation*}
E(Y^2)=\int_0^1\int_0^1E(g(x(t)g(x(s))\mathrm{d}t\mathrm{d}s=E(g(x(t))E(g(x(s))=E^2(Y).
\end{equation*}

(ii) Next we prove the formula of $EY_k$. Without loss of generality, we only discuss the case $k=1$, and suppose  $I=[0,a]$ where $a$ is a rational number or an irrational number. When $a$ is an rational number, that is, $a=\frac{r}{s}, (r,s)=1$, we divide $[0,1]$ by $sn$ equal parts, and hence the discretization of the functional $Y=\int_0^{\frac{r}{s}}g(x(t))\mathrm{d}t$ is given by $Y_n=\frac{1}{sn}\sum_{i=1}^{rn}g(x_i)$. So we have
\begin{equation}
EY_n=\frac{1}{sn}\sum_{i=1}^{rn}Eg(x_i)=\frac{r}{s}E(g(x_1)).
\end{equation}
When $a$ is an irrational number, we can use the rational numbers sequence $\{a_r\}$ to approach it, and then by limitation, we get the result.

Further we have $DY_k=0$ by the same reason with the case (i).

(iii) Similar to case (i), with operations of combinatorics, by the theorem 3, we have $E(Y^k)=E^k(Y)$ for any positive integer $k$. Therefore, we have for any $m\geq1$,
\begin{equation}
\lim_{n\rightarrow+\infty}E(Y_n-EY)^m=0.
\end{equation}
In fact, we have
\begin{equation}
(Y_n-EY)^m=\sum_{j=0}^m\binom{m}{j}(-1)^j(EY_n)^{m-j}(EY)^j,
\end{equation}
and hence from $E(Y^k)=\lim_{n\rightarrow+\infty}E(Y_n^k)=E^k(Y)$,
\begin{equation}
\lim_{n\rightarrow+\infty}E(Y_n-EY)^m=\sum_{j=0}^m\binom{m}{j}(-1)^j(EY)^{m-j}(EY)^j=(1-1)^m=0.
\end{equation}
Next, by definition 5,
\begin{equation}
Eh(Y)=\lim_{n\rightarrow+\infty}Eh(Y_n)=\lim_{n\rightarrow+\infty}Eh(\frac{1}{n}\sum_{k=1}^ng(x_k)).
\end{equation}
Further, from the mean value theorem, there is a point $\xi_n\in(Y_n,EY)$ or $\xi_n\in(EY, Y_n)$  such that
\begin{equation*}
h(Y_n)=h(EY)+h'(EY)(Y_n-EY)+\cdots
\end{equation*}
\begin{equation}
+\frac{h^{(k-1)}(EY)}{(k-1)!}(Y_n-EY)^{k-1}+\frac{h^{(k)}(\xi_n)}{(k)!}(Y_n-EY)^{k},
\end{equation}
and then
\begin{equation*}
Eh(Y_n)=Eh(EY)+h'(EY)E(Y_n-EY)+\cdots
\end{equation*}
\begin{equation}
+\frac{h^{(k-1)}(EY)}{(k-1)!}E(Y_n-EY)^{k-1}+E(\frac{h^{(k)}(\xi_n)}{(k)!}(Y_n-EY)^{k}).
\end{equation}
For the last term, we have as $n\rightarrow+\infty$,
\begin{equation}
|E(\frac{h^{(k)}(\xi_n)}{(k)!}(Y_n-EY)^{k})|<A|E(Y_n-EY)^{k}|\rightarrow0.
\end{equation}
So, taking the limitation as $n$ tending to $+\infty$, and using $\lim_{n\rightarrow+\infty}E(Y_n-EY)^m=0$,
we get $Eh(Y)=h(EY)$. The proof is completed.

In general, we have the following theorem.

\textbf{Theorem 5}. For $n$ functions $g_i$ of $m$ variables, suppose that every function $G_i$ of $g_i$ satisfies the same conditions in theorem 3. Take the ball $M=\{x|\int_0^1x^p(t)\mathrm{d}t\leq R^p\}$ in $C[0,1]$ with norm $||x||_p$ where $p\geq1$ when  $p=\frac{p_0}{q_0}$ where $p_0$ is even and $(p_0,q_0)=1$.  Consider the functionals $Y_1,\cdots, Y_n$ with the form
 \begin{equation}
Y_k=\int_{I_{k1}}\cdots\int_{I_{km_k}}g_k(x(t_1),\cdots,x(t_{m_k}))\mathrm{d}t_1\cdots\mathrm{d}t_{m_k},
\end{equation}
where $I_{kj}$ for $ j=1,\cdots,m_k, k=1,\cdots,n$ are subintervals of  $[0,1]$. We have

(i) The average values satisfy
  \begin{equation}
 EY_k=\frac{vol(I_{k1}\times\cdots\times I_{km_k})}{2^m\Gamma^m(\frac{1}{p})p^{\frac{m}{p}-m}}\int_{-\infty}^{+\infty}\cdots\int_{-\infty}^{+\infty}
g_k(Rx_1,\cdots,Rx_m)\mathrm{e}^{-\frac{x_1^p}{p}-\cdots-\frac{x_m^p}{p}}\mathrm{d}x_1\cdots\mathrm{d}x_m,
\end{equation}
\begin{equation}
DY_k=0,
\end{equation}
 where $vol(I_{k1}\times\cdots\times I_{km_k})=l(I_{k1})\times\cdots\times l(I_{km_k})$.

 (ii). Given a function $h$ of $n$ variables, $h\in C^k(-\infty, +\infty)^n$, and there is a positive constant $A$ such that $|\frac{\partial^k h}{\partial x_j^k}|<A$,  we have the nonlinear exchange formula for the  average value of $h(Y_1,\cdots,Y_n)$ on $M$,
\begin{equation}
Eh(Y_1,\cdots,Y_n)=h(EY_1,\cdots,EY_n).
\end{equation}

\textbf{Proof}. (i). Without loss of generality, we only take $n=2$. The proofs of (i) in the theorem can be give by the theorem 3 and the similar  method in theorem 4. For the general case of $I_k$ being subinterval of $[0,1]$, we only need a proof of the following formula for $Y_1=\int_{I_1}g_1(x(t))\mathrm{d}t$ and $Y_2=\int_{I_2}g_2(x(t))\mathrm{d}t$ where $I_1=[0,\frac{r}{s}]$ and $I_2=[\frac{r_1}{s}, \frac{r_2}{s}]$ being the intervals of $[0,1]$ where $r,r_1,r_2,s$ all are positive integers,
\begin{equation}
E(Y_1Y_2)=E(Y_1)E(Y_2).
\end{equation}
 Indeed, we divide $[0,1]$ into $sn$ equal parts, and then
\begin{equation}
Y_{1n}=\frac{1}{sn}\sum_{i=1}^{rn}g_1(x_i),
\end{equation}
\begin{equation}
Y_{2n}=\frac{1}{sn}\sum_{j=r_1n}^{r_2n}g_2(x_j).
\end{equation}
Therefore,
\begin{equation}
Y_{1n}Y_{2n}=\frac{1}{s^2n^2}\sum_{i=1}^{rn}\sum_{j=r_1n}^{r_2n}g_1(x_i)g_2(x_j).
\end{equation}
Without loss of generality, suppose $r_1<r<r_2$. By symmetry, we have
\begin{equation*}
E(Y_{1n}Y_{2n})=\frac{1}{s^2n^2}\sum_{i=1}^{rn}\sum_{j=r_1n}^{r_2n}E(g_1(x_i)g_2(x_j))
\end{equation*}
\begin{equation}
=\frac{1}{s^2n^2}\{(r(r_1-r_2)n^2+r_2n-1)E(g_1(x_i)g_2(x_j))+(rn-r_2n+1)E(g_1(x_1)g_2(x_1))\}
\end{equation}
By theorem 3 and taking limitation as $n$ approaching to infinity, we get the result. For the case of the endpoints of $I_k$ being irrational numbers, we can prove it by using rational numbers to approach irrational number.

(ii). Here in order to prove the nonlinear exchange formula, we only consider the case of $n=2$ and $I_1=I_2=[0,1]$, while general case can be similarly proven. By definition 5,
\begin{equation}
Eh(Y_1,Y_2)=\lim_{n\rightarrow+\infty}Eh(Y_{1n},Y_{2n})=\lim_{n\rightarrow+\infty}Eh(\frac{1}{n}\sum_{i=1}^ng_1(x_i), \frac{1}{n}\sum_{j=1}^ng_2(x_j)).
\end{equation}
Further, from the mean value theorem, there are $\xi_{1n}$ and $\xi_{jn}$ being respectively between $EY_j$ and $Y_{jn}$ for $j=1,2$, that is $\xi_n\in(Y_n,EY)$ or $\xi_n\in(EY, Y_n)$  such that
\begin{equation*}
h(Y_{1n},h_{2n})=h(EY_1,EY_2)+\frac{\partial h}{\partial Y_1 }(EY_1,EY_2)(Y_{1n}-EY_1)+\frac{\partial h}{\partial Y_2 }(EY_1,EY_2)(Y_{2n}-EY_2)
\end{equation*}
\begin{equation*}
+\cdots+\frac{1}{(k-1)!}\sum_{i=0}^{k-1}\frac{\partial^{k-1} h(EY_1,EY_2)}{\partial Y_1^i\partial Y_2^{k-1-i} }(Y_{1n}-EY_1)^i(Y_{2n}-EY_2)^{k-1-i}
\end{equation*}
\begin{equation}
+\frac{1}{(k)!}\sum_{i=0}^{k}\frac{\partial^{k} h(\xi_{1n},\xi_{2n})}{\partial Y_1^i\partial Y_2^{k-i} }(Y_{1n}-EY_1)^i(Y_{2n}-EY_2)^{k-i}.
\end{equation}
And then
\begin{equation*}
Eh(Y_{1n},h_{2n})=h(EY_1,EY_2)+\frac{\partial h}{\partial Y_1 }(EY_1,EY_2)E(Y_{1n}-EY_1)+\frac{\partial h}{\partial Y_2 }(EY_1,EY_2)E(Y_{2n}-EY_2)
\end{equation*}
\begin{equation*}
+\cdots+\frac{1}{(k-1)!}\sum_{i=0}^{k-1}\frac{\partial^{k-1} h(EY_1,EY_2)}{\partial Y_1^i\partial Y_2^{k-1-i} }E((Y_{1n}-EY_1)^i(Y_{2n}-EY_2)^{k-1-i})
\end{equation*}
\begin{equation}
+\frac{1}{(k)!}\sum_{i=0}^{k}E\{\frac{\partial^{k} h(\xi_{1n},\xi_{2n})}{\partial Y_1^i\partial Y_2^{k-i} }(Y_{1n}-EY_1)^i(Y_{2n}-EY_2)^{k-i}\},
\end{equation}
Firstly, by the similar method of proving $\lim_{n\rightarrow+\infty}E(Y_n-EY)^m=0$, we can easily prove $\lim_{n\rightarrow+\infty}E((Y_{1n}-EY_1)^i(Y_{2n}-EY_2)^j)=0$.  Then, for the last term, we have as $n\rightarrow+\infty$,
\begin{equation*}
|\sum_{i=0}^{k}E\{\frac{\partial^{k} h(\xi_{1n},\xi_{2n})}{\partial Y_1^i\partial Y_2^{k-i} }(Y_{1n}-EY_1)^i(Y_{2n}-EY_2)^{k-i}\}|
\end{equation*}
\begin{equation}
\leq A|\sum_{i=0}^{k}E\{(Y_{1n}-EY_1)^i(Y_{2n}-EY_2)^{k-i}\}|\rightarrow 0.
\end{equation}
Therefore, taking the limitation as $n$ tending to $+\infty$, and using $\lim_{n\rightarrow+\infty}E((Y_{1n}-EY_1)^i(Y_{2n}-EY_2)^j)=0$,
we get the nonlinear exchange formula $Eh(Y_1,Y_2)=h(EY_1,EY_2)$. The proof is completed.

Similarly, we can prove the following theorem and omit its proof.

\textbf{Theorem 6}. Denote for $m$ variables function $g$
\begin{equation*}
G(z_1,\cdots,z_{m-1},z)=g(Rz^{\frac{1}{p}}_1,\cdots,R(z-z_1-\cdots-z_{m-1})^{\frac{1}{p}})
\end{equation*}
\begin{equation}
\times\{z_1(z_2-z_1)\cdots(z-z_1-\cdots-z_{m-1})\}^{\frac{1}{p}-1},
\end{equation}
and
\begin{equation}
f_n(z)=\int_0^n\cdots\int_0^nG(z_1,\cdots,z_{m-1},z))\mathrm{d}z_1\cdots\mathrm{d}z_{m-1}
\end{equation}
\begin{equation}
f(z)=\int_0^{+\infty}\cdots\int_0^{+\infty}G(z_1,\cdots,z_{m-1},z))\mathrm{d}z_1\cdots\mathrm{d}z_{m-1}
\end{equation}
Suppose that $g$ satisfies one of the following two conditions

$H_1$. (Differential condition):  $\lim_{n\rightarrow+\infty}f^{(j)}_n(x)=f^{(j)}(x)$ for $j=0,\cdots,k$, $f_n(x)\in C^k[0,+\infty)$, $\lim_{x\rightarrow+\infty}f_n^{(j)}(x)\mathrm{e}^{-x}=0$ for $j=0,\cdots, k-1$, and there is a positive constant number $A$ such that  $|f_n^{(k)}(x)|\leq A$;

$H_2$. (Integral condition): Suppose that  $\lim_{n\rightarrow+\infty}f_n(x)=f(x)$ uniformly for $x$, and
\begin{equation}
\int_0^{+\infty}|f(x)|\mathrm{e}^{-\frac{x}{p}}\mathrm{d}x<+\infty,
\int_0^{+\infty}x^2|f(x)|\mathrm{e}^{-\frac{x}{p}}\mathrm{d}x<+\infty.
\end{equation}

Take the "first quadrant" $M^+=\{x|x(t)\geq0, \int_0^1x^p(t)\mathrm{d}t\leq R^p\}$ of $M$ where $p\geq1$ or in special $p=\frac{p_0}{q_0}$ where $p_0$ is odd and $(p_0,q_0)=1$. Then the average value and variance on $M$ of the functional
\begin{equation}
Y=\int_{I_1}\cdots \int_{I_m}g(x(t_1),\cdots, x(t_m))\mathrm{d}t_1\cdots\mathrm{d}t_m,
\end{equation}
satisfy
\begin{equation}
EY=\frac{l(I_1)\times\cdots\times l(I_m)}{\Gamma^m(\frac{1}{p})p^{\frac{m}{p}-m}}\int_{0}^{+\infty}\cdots\int_{0}^{+\infty}
g(Rx_1,\cdots,Rx_m)\mathrm{e}^{-\frac{x_1^p}{p}-\cdots-\frac{x_m^p}{p}}\mathrm{d}x_1\cdots\mathrm{d}x_m,
\end{equation}
\begin{equation}
DY=0,
\end{equation}
where $I_i$ for $i=1,\cdots m$ are subintervals of $[0,1]$, $l(I_i)$ is the length of $I_i$.
Further, given a function $h$ of $n$ variables, $h\in C^k(-\infty, +\infty)^n$, and there is a positive constant $A$ such that $|\frac{\partial^k h}{\partial x_j^k}|<A$,   we have the nonlinear exchange formula for the average value of $h(Y_1,\cdots,Y_n)$ on $M$,
\begin{equation}
Eh(Y_1,\cdots,Y_n)=h(EY_1,\cdots,EY_n),
\end{equation}
where $Y_k$ is the same with that in theorem 5.

 Formally, the above theorems 3-6 show the complete concentration of measure phenomenon which means that a functional takes its average value on an infinite-dimensional ball with probability 1. However, we have no measure and hence no corresponding probability. Therefore, instead of measure, we use $DY=0$ to character the so-called concentration of measure phenomenon.

 If we consider infinite-dimensional sphere as the limitation of finite dimensional balls, we have formally the routine concentration of measure. Indeed, for the sphere $M(R)=\{x|\int_0^1x^p(t)\mathrm{d}t\leq R^p\}$ ($p=\frac{p_0}{q_0}$ where $p_0$ is even and $(p_0,q_0)=1$) or $M(R)=\{x|\int_0^1x^p(t)\mathrm{d}t\leq R^p, x(t)\geq0\}$ ($p=\frac{p_0}{q_0}$ where $p_0$ is odd and $(p_0,q_0)=1$) in $C[0,1]$ with $p-$norm, we can easily see that the measure concentrates completely on the surface of sphere $M$ since $\frac{V(M(R))}{V(M(r))}=\lim_{n\rightarrow+\infty}(\frac{r}{R})^n=0$ for $r<R$. Therefore, the average value of the functional in the set $M$ (or $M^+$) can be reduced to and equal to the average value on the surface of $M$(or $M^+$).

 It is easy to see that other functionals don't always satisfy $DY=0$. For example, for the functional
\begin{equation}
Y=f(x)=\int_{I_1}\cdots\int_{I_m}g(x(t_1),\cdots,x(t_m),t_1,\cdots,t_m)\mathrm{d}t_1\cdots\mathrm{d}t_m,
\end{equation}
in general, the average value on $M$ or $M^+$ will be respectively
\begin{equation*}
E(Y)=\frac{l(I_1)\times\cdots\times l(I_m)}{2^m\Gamma^m(\frac{1}{p})p^{\frac{m}{p}-m}}\int_{I_1}\cdots\int_{I_m}\int_{-\infty}^{+\infty}\cdots\int_{-\infty}^{+\infty}
g(Rx_1,\cdots,Rx_m,t_1,\cdots,t_m)
\end{equation*}
\begin{equation}
\times\mathrm{e}^{-\frac{x_1^p}{p}-\cdots-\frac{x_m^p}{p}}
\mathrm{d}x_1\cdots\mathrm{d}x_m\mathrm{d}t_1\cdots\mathrm{d}t_m,
\end{equation}
or
\begin{equation*}
E(Y)=\frac{l(I_1)\times\cdots\times l(I_m)}{\Gamma^m(\frac{1}{p})p^{\frac{m}{p}-m}}\int_{I_1}\cdots\int_{I_m}\int_{0}^{+\infty}\cdots\int_{0}^{+\infty}
g(Rx_1,\cdots,Rx_m,t_1,\cdots,t_m)
\end{equation*}
\begin{equation}
\times\mathrm{e}^{-\frac{x_1^p}{p}-\cdots-\frac{x_m^p}{p}}
\mathrm{d}x_1\cdots\mathrm{d}x_m\mathrm{d}t_1\cdots\mathrm{d}t_m.
\end{equation}
But in general $DY\neq0$ and then the complete concentration of measure don't hold.

\textbf{Remark 3}. How to compute the average value  is dependent to how to define it. From the theorems 1 and 2, we have obtained the densities of the coordinates of the points in $M$ and $M^+$, thus it seems that a natural way of getting the average value of $Y$ is just to define it as the formulas (125) and (126) respectively only if the  integrals in the right sides exist and are finite or infinite. However, from the above theorems 3-6,  we can see that we need some restrictive conditions on $g$ to prove the corresponding formulas of average values for integral form functionals. This is because that we use the same discretization of interval $[0,1]$ of $t$ in both $M$ and $Y$ such that we must deal with simultaneously the  limitation of $n$ appearing in integrand and upper limit of integral, according to the definitions 3 and 4. For the purpose of consistency, this is a necessary step.

\section{Further results and discussions}

In the section, I will give some further results and discussions on average values of functionals for some general cases, including the average value depending on discretiization, the definition of average values of the general continuous functionals, and the definition  of average values on the whole space $C[0,1]$ with $p-$norm. Our discussions will focus on the $M$, while it is similar for $M^+$. In addition, we mainly discuss the case of $g$ is a single variable function, and the general case can be easily given. Specially, we prove an important result which means that the average value depends on the partition of $[0,1]$.

\subsection{The average value depends on the discretization}
From these definitions of average value of functionals, it seems that the average value should be independent to the partition of $[0,1]$. But it is not the case.

It is well known that the stochastic integral $\int_a^bX(t)\mathrm{d}W(t)$ is defined by the mean square limitation of the following summation
\begin{equation}
\sum_{k=1}^nX(t'_k)(W(t_k)-W(t_{k-1})),
\end{equation}
where $a=t_0<t_1<\cdots<t_n=b$, and $t'_k\in[t_{k-1},t_k)$. However, if $t_k'$ is taken arbitrarily, the limitation doesn't exist. Therefore, in general, we take $t'_k=t_{k-1}$ to give the Ito integral, and take $t'_k=\frac{t_k+t_{k-1}}{2}$ to give Stratonovich integral. Being analogous to the stochastic integral, our definition of average values of functionals also depends on the discritizations.

For any partition of $[0,1]$, $0=t_0<t_1<\cdots<t_{n-1}<t_n=1$, we use discrete summation $\sum_{k=1}^ng(x_k)\Delta t_k$ as $\lambda\rightarrow 0$ to replace the integral $Y=f(x)=\int_0^1g(x(t))\mathrm{d}t$ where
$\Delta t_k=t_{k}-t_{k-1}$, $\lambda=\max_{1\leq k\leq n}\Delta t_k$ and $x_k=x(\xi_k)$ for arbitrary number $\xi_k\in(t_{k-1},t_{k})$.This is reasonable by the continuity of $g$. The discretization $M_n$ of the ball $M$ is given by $M_n=\{(x_1,\cdots,x_n)|x_1^p\Delta t_1+\cdots+x_n^p\Delta t_n\leq R^p\}$. We will obtain a theorem and then give an example to show that the different dsicretization yields the different average value.

\textbf{Theorem 7}. Take the ball $M=\{x|\int_0^1x^p(t)\mathrm{d}t\leq R^p\}$ in $C[0,1]$ with norm $||x||_p$ where $p\geq1$ when  $p=\frac{p_0}{q_0}$ where $p_0$ is even and $(p_0,q_0)=1$. Suppose that the partition $0=t_0<t_1<\cdots<t_{n-1}<t_n=1$ of $[0,1]$ satisfies the following conditions: $\Delta t_k=\frac{r_j}{n}$ for $k=s_{j-1}n+1,\cdots,s_jn, j=1,\dots,m$, where $s_0=0, s_j>0, r_j>0$, $s_1+\cdots+s_m=1$ and $s_1r_1+\cdots+s_mr_m=1$. Then the average value of functional $Y=\int_0^1g(x(t))\mathrm{d}t$
on the ball $M$ is given by
\begin{equation}
EY=\frac{1}{2R\Gamma (\frac{1}{p})p^{\frac{1}{p}-1}}\int_{-\infty}^{+\infty}\{\sum_{k=1}^ms_kr_kg(r_k^{-\frac{1}{p}}x)\}
\mathrm{e}^{-\frac{x^p}{pR^p}}\mathrm{d}x,
\end{equation}
where $g$ satisfies the same conditions in theorem 4. This means that the average values of functionals depend on the discretizations.

\textbf{Proof}. The main ideas and steps of the proof are similar to the theorem 4, so we only give main steps without detailed process. By the simultaneous discretization of $Y$ and $M$ as above, we have
 the average value of $Y_n$ as
\begin{equation}
EY_n=\lim_{\lambda\rightarrow 0}\frac{\int_{M_n} \sum_{k=1}^ng(x_k)\mathrm{d}v_n\Delta t_k}{\int_{M_n} \mathrm{d}v_n}
\end{equation}
Since $\lambda\rightarrow 0$ implying $n\rightarrow \infty$, we have
\begin{equation}
EY=\lim_{n\rightarrow \infty}\sum_{k=1}^n\frac{\int_{M_n} g(x_k)\mathrm{d}v_n\Delta t_k}{\int_{M_n} \mathrm{d}v_n}.
\end{equation}
By the similar process with the theorem 4, we get
\begin{equation*}
EY_n=\frac{1}{2R\Gamma (\frac{1}{p})p^{\frac{1}{p}-1}}\sum_{k=1}^m\int_{0}^{R(\Delta t_k)^{-\frac{1}{p}}}
\{g(x_k)+g(-x_k)\}
(1-\frac{x_k^p}{R^p}\Delta t_k)^{\frac{n-1}{p}}(n\Delta t_k)^{\frac{1}{p}}\mathrm{d}x_k\Delta t_k
\end{equation*}
\begin{equation*}
=\frac{1}{2R\Gamma (\frac{1}{p})p^{\frac{1}{p}-1}}\sum_{j=1}^m\int_{0}^{Rn^{\frac{1}{p}}}
\{g(x_k(n\Delta t_k)^{-\frac{1}{p}})+g(-x_k(n\Delta t_k)^{-\frac{1}{p}})\}
(1-\frac{x_k^p}{nR^p})^{\frac{n-1}{p}}\mathrm{d}x_k\Delta t_k
\end{equation*}
\begin{equation}
=\frac{1}{2R\Gamma (\frac{1}{p})p^{\frac{1}{p}-1}}\int_{0}^{Rn^{\frac{1}{p}}}
\sum_{j=1}^m s_kr_k\{g(x_kr_k^{-\frac{1}{p}})+g(-x_kr_k^{-\frac{1}{p}})\}
(1-\frac{x_k^p}{nR^p})^{\frac{n-1}{p}}\mathrm{d}x_k.
\end{equation}
By lemma 2, we take the limitation as $n$ tending to infinity and get the conclusion. The proof is completed.

From the above theorem, if $g(x)\neq \sum_{j=1}^ms_jr_jg(r_j^{-\frac{1}{p}}x)$, we will obtain two different average values of $EY$. For example, we take $m=2, s_1=s_2=\frac{1}{2}, r_1=\frac{1}{2}, r_2=
\frac{3}{2}$, then we get from theorem 7
\begin{equation}
(EY)_1=\frac{1}{2R\Gamma (\frac{1}{p})p^{\frac{1}{p}-1}}\int_{-\infty}^{+\infty}\{\frac{1}{4}g(2^{\frac{1}{p}}x)
+\frac{3}{4}g((\frac{2}{3})^{\frac{1}{p}}x)\}
\mathrm{e}^{-\frac{x^p}{pR^p}}\mathrm{d}x,
\end{equation}
and from theorem 4 or from theorem 7 for $m=1, s_1=r_1=1$,
\begin{equation}
(EY)_2=\frac{1}{2R\Gamma (\frac{1}{p})p^{\frac{1}{p}-1}}\int_{-\infty}^{+\infty}g(x)
\mathrm{e}^{-\frac{x^p}{pR^p}}\mathrm{d}x.
\end{equation}
In general $(EY)_1\neq (EY)_2$.

In the paper, the reason of taking uniform partition is based on the symmetry and simplicity and maybe esthetics.
From the theorem, we also see the complexity of the topic.

\subsection{The average values of analytic functionals}

 Firstly, we need the following theorem. Its proof is similar to the theorem 4, so we omit it for simplicity.

\textbf{Theorem 8}. Suppose that the function $a(t_1,\cdots,t_m)$ is integrable, and $g$ satisfies the same conditions in theorem 5. Take the ball $M=\{x|\int_0^1x^p(t)\mathrm{d}t\leq R^p\}$ in $C[0,1]$ with norm $||x||_p$ where $p\geq1$ when  $p=\frac{p_0}{q_0}$ where $p_0$ is even and $(p_0,q_0)=1$.  Then for the functional
 \begin{equation}
Y=\int_{0}^1\cdots\int_{0}^1a(t_1,\cdots,t_m)g(x(t_1),\cdots,x(t_{m_k}))\mathrm{d}t_1\cdots\mathrm{d}t_{m},
\end{equation}
 we have
  \begin{equation*}
 EY=\frac{1}{2^m\Gamma^m(\frac{1}{p})p^{\frac{m}{p}-m}}\int_{0}^{1}\cdots\int_{0}^{1}
a(t_1,\cdots,t_m)\mathrm{d}t_1\cdots\mathrm{d}t_m
\end{equation*}
 \begin{equation}
\times\int_{-\infty}^{+\infty}\cdots\int_{-\infty}^{+\infty}
g(Rx_1,\cdots,Rx_m)\mathrm{e}^{-\frac{x_1^p}{p}-\cdots-\frac{x_m^p}{p}}\mathrm{d}x_1\cdots\mathrm{d}x_m.
\end{equation}

Next we consider the average values of analytic functionals on the ball $M$. According to [4], the analytic functional $Y=f(x)$ can be expanded as the power series
\begin{equation}
f(x)=\sum_{m=0}^{+\infty}\frac{1}{m!}\int_{0}^{1}\cdots\int_{0}^{1}
\frac{\delta^mf(0)}{\delta x(t_1)\cdots\delta x(t_m)}x(t_1)\cdots x(t_m)\mathrm{d}t_1\cdots\mathrm{d}t_m,
\end{equation}
where $\frac{\delta f}{\delta x(t)}$ is the Fr\'{e}chet derivative. By discretization, we have
\begin{equation}
Y_n=f_n(x_1,\cdots,x_n)=\sum_{m=0}^{+\infty}\frac{1}{m!n^m}\sum_{i_1=1}^{n}\cdots\sum_{i_m=1}^{n}a(t_{i_1},\cdots,t_{i_m})
x(t_{i_1})\cdots x(t_{i_m}),
\end{equation}
where we denote $a(t_{i_1},\cdots,t_{i_m})=\frac{\delta^mf(0)}{\delta x(t_{i_1})\cdots\delta x(t_{i_m})}$. If we suppose that $Y_n$ satisfies some strict conditions such as uniform convergence and uniform continuousness, from the theorem 8, we can get the formula of the average value of $Y$ on the ball $M$
\begin{equation*}
EY=\sum_{m=0}^{+\infty}\frac{1}{m!}\int_{0}^{1}\cdots\int_{0}^{1}
\frac{\delta^mf(0)}{\delta x(t_1)\cdots\delta x(t_m)}\mathrm{d}t_1\cdots\mathrm{d}t_m
\end{equation*}
\begin{equation}
\times\{\frac{1}{2R\Gamma (\frac{1}{p})p^{\frac{1}{p}-1}}\int_{-\infty}^{+\infty}x\mathrm{e}^{-\frac{x^p}{pR^p}}\mathrm{d}x\}^m.
\end{equation}
In fact, if for any fixed $n$, $Y_n=f_n(x_1,\cdots,x_n)$ as the functions series is uniform convergent in $M_n$ on the variables $(x_1,\cdots,x_n)$, we can exchange the order of summation and expectation operations, that is
\begin{equation}
EY_n=\sum_{m=0}^{+\infty}\frac{1}{m!n^m}\sum_{i_1=1}^{n}\cdots\sum_{i_m=1}^{n}a(t_{i_1},\cdots,t_{i_m})
E(x(t_{i_1})\cdots x(t_{i_m})),
\end{equation}
 and then if  $EY_n$ as a constant series is uniform convergent for all $n$, we have
\begin{equation*}
\lim_{n\rightarrow+\infty}EY_n=\sum_{m=0}^{+\infty}\lim_{n\rightarrow+\infty}\frac{1}{m!n^m}\sum_{i_1=1}^{n}\cdots\sum_{i_m=1}^{n}a(t_{i_1},\cdots,t_{i_m})
E(x(t_{i_1})\cdots x(t_{i_m}))
\end{equation*}
\begin{equation*}
=\sum_{m=0}^{+\infty}\frac{1}{m!}\int_{0}^{1}\cdots\int_{0}^{1}
\frac{\delta^mf(0)}{\delta x(t_1)\cdots\delta x(t_m)}\mathrm{d}t_1\cdots\mathrm{d}t_m
\end{equation*}
\begin{equation}
\times\{\frac{1}{2R\Gamma (\frac{1}{p})p^{\frac{1}{p}-1}}\int_{-\infty}^{+\infty}x\mathrm{e}^{-\frac{x^p}{pR^p}}\mathrm{d}x\}^m.
\end{equation}

\subsection{The average values on bounded set and the whole space}

From the average values on the ball $M$, we can compute the average values on the  bounded set $M_0=\{x|r\leq||x||_p\leq R\}$ and define the average values on the whole space $C[0,1]$. Denote $E_MY$ to be the average value of $Y$ on $M$. In fact, we have

\textbf{Theorem 9}. Take the ball $M_0=\{x|r^p\leq\int_0^1x^p(t)\mathrm{d}t\leq R^p\}$ in $C[0,1]$ with norm $||x||_p$ where $p\geq1$ when  $p=\frac{p_0}{q_0}$ where $p_0$ is even and $(p_0,q_0)=1$, and denote $M_R=\{x|\int_0^1x^p(t)\mathrm{d}t\leq R^p\}$ and $M_r=\{x|\int_0^1x^p(t)\mathrm{d}t\leq r^p\}$, $r<R$. Suppose $g$ satisfies the same conditions in theorem 4. Then for the functional
 \begin{equation}
Y=\int_{0}^1g(x(t))\mathrm{d}t,
\end{equation}
we have
\begin{equation}
E_{M_0}Y=E_{M_R}Y.
\end{equation}

\textbf{Proof}. By the same discretization of $x(t)$ on $[0,1]$, we have $Y_n=\frac{1}{n}\sum_{k=1}^n g(x_k)$,
$M_{0n}=\{(x_1,\cdots,x_n)|nr^p\leq x_1^p+\cdots+x_n^p\leq nR^p\}$, $M_{1n}=\{(x_1,\cdots,x_n)| x_1^p+\cdots+x_n^p\leq nr^p\}$ and $M_{1n}=\{(x_1,\cdots,x_n)| x_1^p+\cdots+x_n^p\leq nR^p\}$. Then
\begin{equation*}
E_{M_{0n}}=\frac{\int_{M_{0n}}Y_n\mathrm{d}v_n}{\int_{M_{0n}}\mathrm{d}v_n}
=\frac{\int_{M_{2n}}Y_n\mathrm{d}v_n-\int_{M_{1n}}Y_n\mathrm{d}v_n}
{\int_{M_{2n}}\mathrm{d}v_n-\int_{M_{1n}}\mathrm{d}v_n}
\end{equation*}
\begin{equation}
=\frac{\frac{\int_{M_{2n}}Y_n\mathrm{d}v_n}{\int_{M_{2n}}\mathrm{d}v_n}-\frac{\int_{M_{1n}}\mathrm{d}v_n}{
\int_{M_{2n}}\mathrm{d}v_n}\frac{\int_{M_{1n}}
Y_n\mathrm{d}v_n}{\int_{M_{1n}}\mathrm{d}v_n}}
{1-\frac{\int_{M_{1n}}\mathrm{d}v_n}{\int_{M_{2n}}\mathrm{d}v_n}}.
\end{equation}
From as $n\rightarrow+\infty$,
\begin{equation}
\frac{\int_{M_{1n}}\mathrm{d}v_n}{\int_{M_{2n}}\mathrm{d}v_n}=(\frac{r}{R})^n\rightarrow 0,
\end{equation}
we have
\begin{equation}
E_{M_0}Y=\lim_{n\rightarrow+\infty}E_{M_{0n}}Y_n=E_{M_R}Y.
\end{equation}
The proof is completed.

From the theorem, we can see that the integrations on the ball is the most basic, by which we not only compute the average values on bounded set, but also we can define and compute the integrations on the whole space  $C[0,1]$. Indeed, a natural definition for the average value of functional $Y$ on the whole space $C[0,1]$ is given as
\begin{equation}
E_{C[0,1]}Y=\lim_{R\rightarrow+\infty}E_{M_R}Y.
\end{equation}
In other words, we first find the average value on the ball with radius $R$, and then let $R$ to tend to infinity.

\subsection{The definition of average value of a general continuous functional}

In previous sections, all results are about the functionals with integral forms. Here we will discuss how to define the average values for a general continuous functionals on the ball $M$. From the theorem 7, we have known that the average values depend on the discretizations, therefore, we only consider the uniform dscretization.  We take the uniform partition $0=t_0<t_1<\cdots<t_{n-1}<t_n=1$ of $[0,1]$, and $t_k=\frac{k}{n}$ for $0\leq k \leq n$. For any $x(t)\in C[0,1]$, we use the piecewise linear function $\widetilde{x_n}(t)$ to replace it, where
\begin{equation}
\widetilde{x_n}(t)=\frac{x(t_k)-x(t_{k-1})}{t_k-t_{k-1}}(t-t_{k-1})+x(t_{k-1}), t_{k-1}\leq t\leq t_k, 0\leq k\leq n.
\end{equation}
Then a continuous functional $f(x)$ can be discretized as $Y_n=f_n(x_1,\cdots,x_n)=f(\widetilde{x_n})$. Now we can define the average value of $f$ on the ball $M$ as follows,
\begin{equation}
EY=\lim_{n\rightarrow+\infty}\frac{\int_{M_{n}}f_n(x_1,\cdots,x_n)\mathrm{d}v_n}{\int_{M_{n}}\mathrm{d}v_n},
\end{equation}
where $M_n$ is the discretization of $M$ in terms of the same partition of $[0,1]$.

If we take the limitation as the radius $R$ of $M$ tending to infinity, we can define the average value of $f$ on the whole space. Although we can give the definition, we cannot give the general method to compute exactly the average values for general continuous functionals. For the concrete functional, we can always find method to compute its average value exactly or approximately.

For those functionals of $x(t)$ and $x'(t)$, such as $Y=f(x,x', \cdots)$, by discretization and replacing derivative $x'$ by difference quotient, we can define and compute their average value on $M$ and $C[0,1]$.

\section{Conclusion}

Essentially, the concept of average value depends on how to define it. Without measure, we also give such definitions on infinite dimensional balls. In particular, we prove that the average value depends on the discretization. This is an important and interesting result from which we know that the definition and computation of average values of functional on $C[0,1]$ with $p-$norm is subtle. In addition, we use variance to describe the deviation of functional from its average value. Borrowing the language of probability, we give versions of geometry for main theorems, and we can also see that it is just the geometrical versions to provide more convenient and simple viewpoint for understanding the infinite dimensional problem. Among those, we use the variance being zero to character concentration without measure. However, a puzzled problem for me is how to give an exact explanation for $DY=0$.  The simplicity of the result is unexpected but reasonable, and hence it is interesting for us. It is well known that in general few dimension means simple. However, from the above results, we can see that in some degree infinite-dimensional integrals are more easy to compute than the finite dimensional integrals. This also means that sometimes infinite dimension contains simplicity. In other words, more free and more simple! This is also the essence of statistical physics. On the other hand, since the average values of functionals depend on the discretizations, there exist subtle complexity in infinite dimensional cases.

\textbf{Acknowledgments}: Thanks to Professor Erich Novak for his helpful comments on theorems 1 and 2 in the first version so that I can improve the whole manuscript. In particular, due to his valuable comments to each version, I can think about the whole topic and then find more results.


\begin{thebibliography}{2}
\bibitem{wiiw}Feynman R P, Hibbs A R.  Quantum mechanics and path integrals: Emended edition. Courier Dover Publications, 2012.
\bibitem{t4}Kleinert H. Path integrals in quantum mechanics, statistics, polymer physics, and financial markets. World Scientific, 2009.
\bibitem{f4}Baaquie B E. Quantum finance: Path integrals and Hamiltonians for options and interest rates. Cambridge University Press, 2007.
  \bibitem{pa}L\'{e}vy P.  Problems concrets d'analyse fonctionnelle. Gauthier-Villars. Paris, 1951.
 \bibitem{wie}Wiener N.  Differential space. Journal of Mathematical Physics. 1923,\textbf{2}: 131-174.
 \bibitem{ew}Wiener N.  The mean of a functional of arbitrary elements. Annals of Mathematics. 1920, \textbf{ 22}: 66-72.
 \bibitem{we1}Wiener N. The average of an analytic functional and the Brownian
 movement. Proceedings of the National Academy of Sciences of the United States of American. 1921, \textbf{7}: 294-298.
 \bibitem{we12}Wiener N. The average of an analytic functional. Proceedings of
  the National Academy of Sciences of the United States of American. 1921, \textbf{7}: 253-260.
 \bibitem{k3}Kac M.  On the average of a certain Wiener functional and a related limit
 theorem in calculus of probability. Transactions of the American Mathematical Society. 1946,\textbf{ 59}: 401-414.
 \bibitem{cm1}Cameron R H, Martin W T. Transformations of Weiner integrals under translations. Annals of Mathematics. 1944,\textbf{45}: 386-396.
   \bibitem{cm2} Cameron R H, Martin W T.  Evaluation of various Wiener integrals by use of certain Sturm-Liouville differential equations. Bulletin of the American Mathematical Society. 1945, \textbf{51}: 73-90.
     \bibitem{cm3} Cameron R H, Martin W T.  Transformations of Wiener integrals under a general class of linear transformations. Transactions of the American Mathematical Society. 1945, \textbf{58}: 184-219.


 \bibitem{g}Gelfand I M, Vilenkin N Y.  Generalized functions. Vol. 4, Applications of harmonic analysis. Moscow: Academic Press. 1964.

 \bibitem{mmm}Albeverio S, Hoegh-Krohn R, Mazzucchi S.  Mathematical theory of Feynman
 path integrals: an introduction. New York: Springer-verlag, 2008.
\bibitem{gy}Skorokhod A V. (1974). Integration in Hilbert space. New York: Springer-Verlag, 1974.
\bibitem{g}Hui-Hsiung Kuo, Gaussian measures in Banach spaces,
Lect. Notes in Math 463, Springer, Berlin, 1975.
\bibitem{er}Giles M B. Multilevel Monte Carlo path simulation, Oper. Res. 2008, 56: 607-617.
\bibitem{e2}Giles M B, Waterhouse B J. Multilevel quasi-Monte Carlo path simulation. Advanced Financial Modelling, Radon Series on Computational and Applied Mathematics, 2009, 8: 165-181.

\bibitem{n1}Novak E, Wo?niakowski H. Tractability of Multivariate Problems: Standard information for functionals. European Mathematical Society, 2008.
\bibitem{n2}Novak E, Ritter K. High dimensional integration of smooth functions over cubes. Numerische Mathematik, 1996, 75(1): 79-97.

\bibitem{n3}Novak E, Wozniakowski H. When are integration and discrepancy tractable?. LONDON MATHEMATICAL SOCIETY LECTURE NOTE SERIES, 2001: 211-266.

\bibitem{n4}Novak E, Ritter K, Schmitt R, Steinbauer A. On an interpolatory method for high dimensional integration. Journal of computational and applied mathematics, 1999, 112(1-2): 215-228.

\bibitem{e3}Gnewuch  M.  Weighted geometric discrepancies and numerical integration on reproducing
kernel Hilbert spaces, J. Complexity. 2012, 28: 2-17.
\bibitem{e4} Gnewuch M.  Infinite-dimensional integration on weighted Hilbert spaces, Math.
Comp. 2012, 81: 2175-2205.


\bibitem{e6}Wasilkowski G W, W\"{o}zniakowski H. On tractability of path integration, J. Math.
Physics. 1996, 37:2071-2088.
\bibitem{g1}Kuo F Y,  Nuyensb D, Plaskotac L, Sloana I H, Wasilkowski G W.  Infinite-dimensional integration and the multivariate decomposition method. Journal of Computational and Applied Mathematics. 2017,
 \textbf{326}(15):217-234.
 \bibitem{3}Gnewuch M, Mayer S, Ritter K.  On weighted Hilbert spaces and integration of functions of infinitely many variables. Journal of Complexity.   2014), \textbf{30}(2): 29-47.
 \bibitem{g2}Plaskota L, Wasilkowski G W. Efficient algorithms for multivariate and $\infty-$variate integration with exponential weight. Numerical Algorithms.  2014, \textbf{67}(2): 385-403.
\bibitem{4}Dick J, Gnewuch M.  Infinite-dimensional integration in weighted Hilbert spaces: anchored decompositions, optimal deterministic algorithms, and higher-order convergence. Foundations of Computational Mathematics. 2014, \textbf{14}(5): 1027-1077.

\bibitem{E}Milman V D, Schechtman G. Asymptotic Theory of Finite Dimensional
Normed Spaces: Isoperimetric Inequalities in Riemannian Manifolds. New York: Springer-verlag, 1986.
\bibitem{Si}Ledoux M.  The concentration of measure phenomenon.
American Mathematical Soc., 2005.
\bibitem{Li}Talagrand M.A new look at independence. The Annals of probability. 1996,  \textbf{24}: 1-34.

 \bibitem{ff}Fikhtengol'ts G M.  A Course of Differential and Integral Calculus.   Moscow: Science Pres.
1969.

\bibitem{li}Liu Cheng-shi. Maximal non-symmetric entropy leads naturally to zipf's Law. Fractals, 2008, 16(01): 99-101.
\bibitem{li2}Liu Cheng-shi. Nonsymmetric entropy and maximum nonsymmetric entropy principle. Chaos, Solitons and Fractals, 2009, 40(5): 2469-2474.
\bibitem{jl}Liu Cheng-shi. Lectures on the mean values of functionals--An elementary introduction to infinite-dimensional probability. arXiv preprint arXiv:1705.03584, 2017.
\end{thebibliography}
\end{document}